\documentclass[reqno,11pt]{amsart}    

\usepackage{amscd,amsfonts,mathrsfs,amsthm,enumerate}
\usepackage{amssymb, amsmath}

\usepackage{color}           
\definecolor{myowncolor1}{rgb}{0,0,.75}
\definecolor{darkred}{rgb}{0.7,0,0}
\definecolor{darkgreen}{rgb}{0,0.5,0}
\definecolor{darkblue}{rgb}{0,0,.75}
\definecolor{orange}{rgb}{1,.5,0}
\definecolor{papier}{rgb}{1,1,.97}

\ifx\pdfoutput\undefined     
\usepackage{graphicx}
\else \usepackage[pdftex]{graphicx}
\usepackage[pdftex,
            colorlinks, 
            linkcolor=darkred, 
            citecolor=darkgreen, 
            urlcolor=blue]{hyperref}

\fi

\usepackage{html}            
\usepackage{ifthen}          
\usepackage{calc}            
  
\def\nur      {\langle\nu,e_r\rangle}
\def\real     #1{{\mathbb R^{#1}}}
\def\complex  #1{{\mathbb C^{#1}}}

\def\integer  #1{{\mathbb Z^{#1}}}

\def\d        #1#2{{#1}_{#2}}
\def\dd       #1#2#3{{#1}_{#2#3}}
\def\ddd      #1#2#3#4{{#1}_{#2#3#4}}

\def\uu       #1#2#3{{#1}^{#2#3}}

\def\lap      {\Delta }
\def\dt       {\frac{d}{dt}\,}
\def\sE       {\mathscr{E}}

\def\picture#1#2{
\begin{latexonly}
\ifx\pdfoutput\undefined 
    \includegraphics[width=#2\hsize]
      {#1.eps}
\else
    \includegraphics[width=#2\hsize]
      {#1.jpg}
\fi
\end{latexonly}
}

\newtheorem{theorem}{Theorem}[section]   
\newtheorem{lemma}[theorem]{Lemma}   
\newtheorem{corollary}[theorem]{Corollary}   
\newtheorem{proposition}[theorem]{Proposition}   
\newtheorem{remark}[theorem]{Remark}   
  
\theoremstyle{definition}   
\newtheorem{definition}[theorem]{Definition}

\numberwithin{equation}{section}

\newcommand{\bfig}{\begin{figure}}
\newcommand{\efig}{\end{figure}}

\makeatletter
\def\pproof #1{\@ifnextchar[\opargproof
{\opargproof[\it Proof of Theorem {#1}.]}}
\def\opargproof[#1]{\par\noindent {\bf #1 }}
\def\endpproof{{\unskip\nobreak\hfil\penalty50\hskip8mm\hbox{}
\nobreak\hfil
\(\square\)\parfillskip=0mm \par\vspace{3mm}}}
\makeatother

\makeatletter
\def\peroof #1{\@ifnextchar[\opargproof
{\opargproof[\it Proof of Proposition {#1}.]}}
\def\opargproof[#1]{\par\noindent {\bf #1 }}
\def\endperoof{{\unskip\nobreak\hfil\penalty50\hskip8mm\hbox{}
\nobreak\hfil
\(\square\)\parfillskip=0mm \par\vspace{3mm}}}
\makeatother

\newcommand{\cf}{cf.\ }

\newcommand{\set}[1]{\{#1\}}
\newcommand{\setlr}[1]{\left\{#1\right\}} 
\newcommand{\abs}[1]{|#1|}
\newcommand{\lra}{\longrightarrow}

\newcommand{\mtq}[1]{\qquad\text{\rm #1}\;\;}
\newcommand{\mtqq}[1]{\qquad\text{\rm #1}\qquad}

\newcommand{\C}{\mathbb C}
\newcommand{\D}{\mathbb D}

\newcommand{\T}{\mathbb T}
\newcommand{\Z}{\mathbb Z}

\DeclareMathOperator{\rot}{\mathtt{rot}}
\DeclareMathOperator{\wind}{\mathtt{wind}}

\DeclareMathOperator{\mi}{\mu}
\DeclareMathOperator{\pt}{\mathrm{int}}

\DeclareMathOperator{\su}{\mathrm{SU}}
\DeclareMathOperator{\uni}{\mathrm{U}}
\begin{document}
\title
[Mean curvature flow of Lagrangian submanifolds]
{Mean curvature flow of monotone Lagrangian submanifolds}

\author[Groh]{\sc K. Groh}

\author[Schwarz]{\sc M. Schwarz}
\address{Universit\"at Leipzig, Institut f\"ur Mathematik, 
  Augustusplatz 10--11, 04109 Leipzig, Germany}
\email{Matthias.Schwarz@mathematik.uni-leipzig.de} 
\email{Kai.Zehmisch@mathematik.uni-leipzig.de} 

\author[Smoczyk]{\sc K. Smoczyk}
\address{Universit\"at Hannover, Institut f\"ur Differentialgeometrie,    
  Welfengarten 1, 30167 Hannover, Germany}
\email{groh@math.uni-hannover.de} 
\email{smoczyk@math.uni-hannover.de}

\author[Zehmisch]{\sc K. Zehmisch}

\begin{abstract}
We use holomorphic disks to describe the formation of singularities 
in the mean curvature flow of monotone Lagrangian submanifolds in 
$\complex{n}$.
\end{abstract}

\renewcommand{\subjclassname}{%
  \textup{2000} Mathematics Subject Classification}  
   
\date{June 18, 2006.}   
\thanks{Supported by DFG, priority program SPP 1154, SM 78/1-1}
\maketitle

\section{Introduction}
Let $F_0\colon L\to\complex{n}$ 
be a smooth Lagrangian immersion of a compact, oriented manifold $L$ 
without boundary, i.e. $\text{dim}(L)=n$ and $F^*_0\omega=0$, 
where $\omega(\cdot,\cdot)=\langle J\cdot,\cdot\rangle$ is the standard
symplectic form on $\complex{n}$ given by the composition of the complex
structure $J$ and the euclidean metric $\langle\cdot,\cdot\rangle$.

Suppose $[0,T_{\text{sing}})$ is the maximal time interval such that a smooth solution
$F:L\times[0,T_{\text{sing}})\to\complex{n}$ of the mean curvature flow equation
\begin{gather}\label{mcf}
\dt F=\overrightarrow H\tag{MCF}\\
{F}(\cdot,0)=F_0\notag
\end{gather}
exists. It is well known that $0<T_{\text{sing}}<\infty$ and that 
$L_t:=F_t(L)$ is Lagrangian, where we set $F_t:=F(\cdot,t)$. 
Moreover, we have
$$\limsup_{t\to T_{\text{sing}}}|A|^2=\infty,$$
where $A=\nabla dF$ is the second fundamental form.

At the singular time $T_{\text{sing}}$ a singularity will form.
The main aim of this paper is to understand the nature of singularities that monotone Lagrangian
submanifolds might develop under various additional assumptions. This is motivated by a strong
link between monotone Lagrangians and self-similar solutions of the Lagrangian mean curvature 
flow. To explain this relation, let us first recall the notion of blow-up points for the mean 
curvature flow and the notion of monotone Lagrangian submanifolds in symplectic geometry.

To understand what happens as a singularity forms under the mean curvature flow
in euclidean space, one introduces the notion of {\em blow-up points}. 
By definition, a blow-up point $p\in \complex{n}$ 
is a point such that one can find $x\in L$ with $F(x,t)\to p$ as $t\to T_{\text{sing}}$
and $|A|(x,t)$ becomes unbounded as $t\to T_{\text{sing}}$. From the evolution equation of
the second fundamental form one immediately deduces for $T_{\text{sing}}<\infty$ that 
$$\sup_{L_t}|A|^2\ge \frac{c}{T_{\text{sing}}-t}$$
for a positive constant $c$ and all $0\le t< T_{\text{sing}}$.
{F}ollowing Huisken \cite{huisken 1}, a singularity forms under the
so-called {\em type-1 condition}, if
\begin{equation}\label{type-1}\sup_{L_t}|A|^2\le\frac{c}{T_{\text{sing}}-t}
\end{equation}
for some (different) constant $c>0$.
All remaining cases of singularities with a different blow-up rate
are called of {\sl type-2}.

In \cite{huisken 1}, in the analogous context of mean curvature flow
for hypersurfaces in $\real{n+1}$, Huisken has given an answer to the 
question what happens if a singularity forms under the type-1 hypothesis. 
This holds true as well in higher codimension, in particular also for 
Lagrangian submanifolds.

In this situation, a singularity forming in euclidean space looks like
a {\em self-similarly contracting solution} (see definition below) 
after an appropriate rescaling procedure. 
To make this precise in our context, we assume for simplicity 
that the origin $O\in\complex{n}$ is a blow-up point for our flow of Lagrangian
immersions (since the flow is invariant under isometries, we may assume
w.l.o.g. that the blow-up point is given by the origin). 
Then, as in \cite{huisken 1}, we define the rescaled immersions 
$\widetilde F(x,s)$ by
$$\widetilde F(x,s):=(2(T_{\text{sing}}-t))^{-1/2}F(x,t),$$
$$s(t):=-\frac{1}{2}\log (T_{\text{sing}}-t).$$
The submanifolds $\widetilde L_s:=\widetilde F(L,s)$ are defined for 
$-\frac{1}{2}\log T_{\text{sing}}\le s<\infty$ and satisfy the equation
\begin{equation}\label{rmcf}
\frac{d}{ds}\widetilde F(x,s)=\widetilde{\overrightarrow H}(x,s)
+\widetilde F(x,s),
\end{equation}
where $\widetilde{\overrightarrow H}$ 
is the mean curvature vector of $\widetilde L_s$. 

Under the type-1 assumption (\ref{type-1}) one can prove as in \cite{huisken 1} 
that for each sequence
$s_j\to\infty$ there is a subsequence $s_{j_k}$ such that $\widetilde
L_{s_{j_k}}$ converges to an immersed non-empty smooth limiting
Lagrangian $\widetilde L_\infty$. In particular, this limit Lagrangian is a 
self-similarly contracting solution. 
A solution $F(x,t)$ of (\ref{mcf}) is called {\sl self-similar}, if we can
find a rescaling $\widetilde F(x,t)=\psi(t)F(x,t)$, $\psi(t)\in\real{}$
such that all normal components of $\frac{d}{dt}\widetilde F$ vanish.
This means that the images of $F$ evolve just by homotheties.
Self-similar solutions of the mean curvature flow are determined 
by solutions of the following quasilinear elliptic system
\begin{equation}\label{selfsim}
\overrightarrow H=-\frac{\varepsilon}{2} F\sp\perp,
\end{equation}
with some constant $\varepsilon$. Under the flow, $\varepsilon$ will then vary with time but will
keep its sign. In principle, one has to destinguish 3 cases. A  
solution is called {\sl self-contracting, trivial (i.e. minimal)} or 
{\sl self-expanding} according to whether $\varepsilon >0, \varepsilon=0$ or $\varepsilon<0$.

In the Lagrangian setting, we can reformulate equation (\ref{selfsim}) on $L$.
A self-similar solution must satisfy the following relation between 1-forms 
on $L$
$$
H:=F^*(\omega(\cdot,\overrightarrow H))
=-\frac{\varepsilon}{2}F^*\left(\omega(\cdot,F)\right).
$$
$H$ is the mean curvature form. For the exact symplectic form
$\omega=d\lambda$, $\lambda=\frac{1}{2}\left(xdy-ydx\right)$ on $\complex{n}$, it follows that
\begin{equation}\label{monotone}
  H=\varepsilon F^*\lambda.
\end{equation}
This is due to the fact that rescaling is performed by the flow of
the Liouville vector field $F$, i.e. $\mathcal{L}_F\omega=\omega$.
Note also that the Lagrangian submanifolds under consideration may 
be noncompact.

The mean curvature form $H$ is closed. Moreover, by a result of Morvan
\cite{demaslov} (see also \cite{cieliebak goldstein} for a recent
extension), the complex volume form $dz$ on $\complex{n}$ satisfies
$$dz_{|L}=e\sp{i\alpha}d\mu\,,\quad d\alpha=H,$$
where $\alpha$ is the locally defined {\sl Lagrangian angle} of $L$ and 
$d\mu$ denotes the induced volume form on $L$.

In the Lagrangian context, there is a strong connection between self-similarly
contracting solutions and {\sl monotone} Lagrangian submanifolds.
We recall the following definition (cf. \cite{oh}):

\begin{definition}
Assume $F:L\to\complex{n}$ is Lagrangian. Suppose $u:\D\to\complex{n}$
is a smooth map of the $2$-disk $\D$ whose boundary $\partial u:=u(\partial \D)$
is contained in $L_0:=F(L)$. Let $\mu(u)$ be the Maslov index of $u$ and 
$\omega(u):=\int_{\D}u\sp*\omega$ its symplectic area.
$L_0$ is called {\sl monotone},
if there exists a positive 
constant $\varepsilon>0$ (called {\sl monotonicity constant})
such that
\begin{equation}\label{mon}
\mu(u)=\frac{\varepsilon}{\pi}\,\omega(u)
\end{equation} 
on all such disks $u$.
\end{definition}

\begin{remark}\label{rem 1}
Assume $u$ is a disk with boundary curve $\gamma\subset L$ and that
$F:L\to(\complex{n},\omega=d\lambda)$ is Lagrangian. By Stokes' theorem the
symplectic area is 
$$\omega(u)=\int_\gamma F\sp*\lambda$$
and by Morvan's result \cite{demaslov} we have $H(\gamma)=\pi\mu(u)$ so that
(\ref{mon}) becomes
\begin{equation}\label{monotonicity}
[H]=\varepsilon[F^*\lambda].
\end{equation}
In particular, by (\ref{monotonicity}), self-similarly contracting solutions are monotone.
Moreover, if $L_0$ is monotone with monotonicity constant $\varepsilon$, 
then for $\tilde L_0:=cL_0$ we have $\tilde H=H$, 
$\tilde F^*\lambda=c^2F^*\lambda$, 
so that $\tilde L_0$ is monotone with monotonicity constant 
$\tilde \varepsilon=c^{-2}\varepsilon$. Thus, by a simple rescaling argument
we can always assume that the monotonicity constant $\varepsilon$ of a monotone Lagrangian 
submanifold satisfies $\varepsilon=1$.
\end{remark}

\begin{remark}\label{rem 2}
If $\gamma_1, \gamma_2$ are two closed curves on $L$ and $\gamma_1$ is
homotopic to $\gamma_2$, then $H(\gamma_1)=H(\gamma_2)$ and by Stokes' theorem
and the Lagrangian property of $L$ also 
$\lambda_{|L}(\gamma_1)=\lambda_{|L}(\gamma_2)$. In particular, if
$u:\D\to\complex{n},\, \partial u\subset L$ is a holomorphic disk with boundary
on $L$, then $\lambda_{|L}(\partial u)=\omega(u)>0$ and 
$\lambda_{|L}(\gamma)>0$ for any curve $\gamma$ homotopic to 
$\partial u$. 
\end{remark}

Since self-similarly contracting solutions are monotone and self-similar solutions appear
as type-1 blow-ups, it is natural to ask the following: 

{\sl Does the flow preserve the monotonicity of Lagrangian submanifolds and
will monotone Lagrangians develop type-1 singularities?}

We will see, as stated in the next proposition, 
that the first part of this question has a positive answer 
whereas the second part is wrong in general.

\begin{proposition}\label{lemma mon}
Suppose $F:L\times[0,T_{\text{sing}})\to\complex{n}$ 
is a solution of the Lagrangian mean curvature flow. If
$L_0=F(L,0)$ is monotone with monotonicity constant $\varepsilon_0$, then
$[F^*\lambda]=\frac{1}{\varepsilon_0}(1-\varepsilon_0t)[H]$ for any
$t\in[0,T_{\text{sing}})$.
In particular, $L_t=F(L,t)$ remains monotone for $t\in[0,\frac{1}{\varepsilon_0})$, 
but with the monotonicity constant 
$\varepsilon_t:=\frac{\varepsilon_0}{1-\varepsilon_0t}>0$.
\end{proposition}

{F}rom this proposition and Remark \ref{rem 2} we conclude that 
the areas of holomorphic disks with boundary 
on a monotone Lagrangian torus (if they exist)
all shrink with the same rate. This means
that they would all disappear at the same time $t=\frac{1}{\varepsilon}$. 
Therefore, one is tempted
to conjecture that monotone tori will completely shrink to a single
point and will always develop a type-1 singularity at 
$T_{\text{sing}}=\frac{1}{\varepsilon}$, 
which after blow-up looks like a self-similarly shrinking Lagrangian torus. 
That this is not the case, even under strong additional geometric assumptions,
is somehow surprising. Even under the assumption that the Lagrangian 
submanifolds admit quasi-fillings (see definition below)
for some homology class 
$e\in H_1(L;\integer{})$ one still might encounter type-2
singularities. This will be demonstrated for examples of embedded,
equivariant, monotone $2$-tori that we discuss in this article.

Since the existence of holomorphic disks with boundary on the Lagrangian
submanifold is important, the next theorem due to Gromov \cite{gro} 
for embedded Lagrangian submanifolds is of particular interest.

\begin{theorem}\label{gromov}(Gromov)
Suppose $L\subset\complex{n}$ is an embedded closed Lagrangian. Then
there exists at least one holomorphic disk $u:\D\to\complex{n}$ whose boundary is contained in $L$.
In particular, $[\lambda_{|L}]\neq 0$.
\end{theorem}

In view of Gromov's theorem and of Remark \ref{rem 2} one would like 
to control the embeddedness of Lagrangian submanifolds under the mean 
curvature flow. It is well-known that embeddedness is preserved for 
hypersurfaces moving by their mean curvature and that this, in general,
is no longer true in higher codimensions.
Therefore we define a constant $T_{\text{emb}}$ as follows:

Suppose $F_0:L\to\complex{n}$ is a Lagrangian embedding and that $F:L\times[0,T_{\text{sing}})\to\complex{n}$
is the maximal smooth solution of the mean curvature flow. Then we set
$$T_{\text{emb}}:=\sup\{t\in[0,T_{\text{sing}}):L_\tau \text{ is embedded for } \tau\in[0,t]\}.$$

By definition one has
$$0<T_{\text{emb}}\le T_{\text{sing}}.$$

\def\optional{
One of the aims of this article is to clarify the relation between Hamiltonian 
isotopies of monotone Lagrangian tori in $\complex{n}$
and the behaviour of holomorphic disks in $\complex{n}$ in order to find a
way to control the Lagrangian mean curvature flow.
It turns out, as we will see, that equivariant tori in $\complex{2}$ form a rich source of 
interesting and non-trivial examples in the four-dimensional case.

In general, one may ask why symplectic geometers always talk about Lagrangian 
tori if they think four-dimensionally.
One reason of course is that the two-torus $\T^2$ is a very simple manifold.
But more important is that the only closed orientable surface $F$ which admits 
a Lagrangian embedding into $\complex{2}$ is the torus, because the homological 
self-intersection number $F\cdot F=-\chi(F)$ vanishes (\cf \cite{lastchapter}).

The main questions about Lagrangian tori are of global nature.
How does a  Lagrangian torus sit in $\complex{2}$?
Are there numerical invariants?
How can one distinguish certain types of tori in $\complex{2}$?
We remark that Eliashberg studied compact Lagrangian cylinders 
\cite{topknot} and
Eliashberg and Polterovich studied Lagrangian planes \cite{trivial} in $\complex{2}$.
Their results essentially tell us that there are no local methods to 
answer these questions.

A similar issue arises when one considers Lagrangian immersions instead of 
embeddings.
{F}rom the symplectic point of view, Lagrangian immersions are of less interest
because of Gromov's $h$-principle \cite{partial}.
Roughly it says that all questions about Lagrangian immersions can be 
answered with purely topological methods, see \cite{lastchapter}.
But it is of central importance, for example in view of the Lagrangian 
mean curvature flow,
to decide whether a given regular homotopy of Lagrangian tori
starting with an embedding will be a Lagrangian isotopy.
}

It turns out that embeddedness is important in the context of monotone 
Lagrangian submanifolds. We will prove the following result:

\begin{theorem}\label{theo main 1}
Suppose $F_0:L\to\complex{n}$ is an embedding of a compact, oriented, monotone Lagrangian 
submanifold with monotonicity constant $\varepsilon>0$. 
Then $T_{\text{emb}}\le\frac{1}{\varepsilon}.$
\end{theorem}

The proof of this theorem will be carried out in Section \ref{sec gro} and 
is based on Gromov's result mentioned above. 

\begin{figure}[h]
\begin{center}
\picture{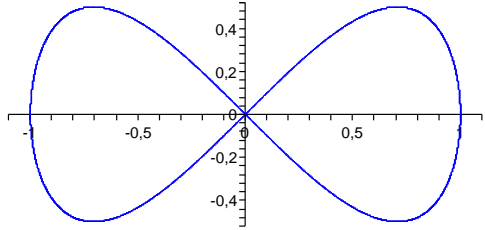}{{.5}}
\caption{\sl For the figure eight curve $\gamma$ 
we have $[\lambda_{|\gamma}]=[H]=0$, so that 
$[H]=\varepsilon[\lambda_{|L}]$ holds for any 
$\varepsilon\in\real{}$. Consequently, $\gamma$ is monotone and the
monotonicity constant is not related to the singular time $T_{\text{sing}}$.}
\label{figure eight}
\end{center}
\end{figure}

A corollary of Theorem \ref{theo main 1} is:
\begin{corollary}\label{cor main 1}
Suppose $F_0:L\to\complex{n}$ is an embedding of a compact, oriented, monotone Lagrangian 
submanifold with monotonicity constant $\varepsilon>0$. 
If $T_{\text{emb}}=T_{\text{sing}}$, 
i.e. if the solution stays embedded for all 
$t<T_{\text{sing}}$, then $T_{\text{sing}}\le\frac{1}{\varepsilon}$.
\end{corollary}

An obvious question is, if
$T_{\text{sing}}>\frac{1}{\varepsilon}$ is possible 
at all, or if for any monotone Lagrangian we must have 
$T_{\text{sing}}\le\frac{1}{\varepsilon}$?
The figure eight curve depicted in Figure \ref{figure eight}
shows that, at least in the immersed case, $\varepsilon$ need not be 
related to $T_{\text{sing}}$.

Another obvious question arising from Corollary \ref{cor main 1}
is, if - under the same assumptions as in Corollary \ref{cor main 1} -
we must have 
$T_{\text{sing}}=\frac{1}{\varepsilon}$ or if 
$T_{\text{sing}}<\frac{1}{\varepsilon}$ is possible.
In this article we will give examples with 
$T_{\text{emb}}=T_{\text{sing}}<\frac{1}{\varepsilon}$ 
that develop type-2 singularities.

On the other hand, it is possible to characterize the solutions with $T_{\text{sing}}=\frac{1}{\varepsilon}$ by the
following proposition (for a proof see Section \ref{sec gro}):

\begin{proposition}\label{prop main 1}
Suppose $F_0:L\to\complex{n}$ is an immersion of a compact, oriented, 
monotone Lagrangian 
submanifold with monotonicity constant $\varepsilon>0$. 
Then $T_{\text{sing}}=\frac{1}{\varepsilon}$, if and only if
the rescaled flow (\ref{rmcf}) is Hamiltonian.
\end{proposition}

In \cite{neves} Neves proved that a Lagrangian submanifold with trivial Maslov
class does not develop a type-1 singularity in finite time. His proof is
based on methods from geometric measure theory. We will give a very elementary
proof for this, in particular we will prove the following extension of 
Neves' result:

\begin{theorem}\label{noself}
Let $F:L\to\complex{n}$ be a Lagrangian self-similar solution, i.e.
a solution of 
$$\overrightarrow H=-\frac{\varepsilon}{2}F^\perp$$
with $\varepsilon\neq 0$
\footnote{$\varepsilon$ can be either positive (contracting solution) or
negative (expanding solution)}. Suppose $L$ has trivial Maslov class
$\mu=\frac{1}{\pi}[H]=\frac{1}{\pi}[d\alpha]=0$ and that there exist 
constants $c>0$,  
$\sigma<\varepsilon$ such that the Lagrangian angle $\alpha$
satisfies $|\nabla^k\alpha|^2\le ce^{\frac{\sigma|F|^2}{4}}$ for 
$0\le k\le 2$. Then $L$ is part of a minimal cone.
\end{theorem}

The proof of this theorem is very short and merely uses partial integration
w.r.t. the Gau\ss\ kernel. Note, that the assumption on $|\nabla^k\alpha|^2$
is automatically satisfied for a type-1 blow-up. Note also that it is an open 
problem in the closed case to determine which Lagrangians admit $\mu=0$.

{F}or the formulation of the next theorem we need the following definition.

\begin{definition}
A type-1 singularity will be called {\sl embedded}, if 
$T_{\text{emb}}=T_{\text{sing}}$
and if in addition, one of the type-1 blow-up limiting Lagrangians 
obtained by the rescaled flow (\ref{rmcf}) is embedded
as well.
\end{definition}

We will prove:

\begin{theorem}\label{theo main 2}
Let $F_0:L\to\complex{n}$ be an embedding of a 
compact, oriented, monotone Lagrangian submanifold with 
monotonicity constant $\varepsilon$. Suppose the Lagrangian submanifold
develops an embedded type-1 singularity. Then 
$T_{\text{sing}}=\frac{1}{\varepsilon}$.
\end{theorem}

However, even under the assumptions made in Theorem \ref{theo main 2}
it is still unclear, if the singular set consists of a single point and
if the blow-up converges to a compact self-similar Lagrangian. This
very much depends on the "families" of holomorphic disks that can be attached
to the evolving submanifolds. For example, in the easiest case of a Lagrangian
$L\subset\complex{2}$ it is in general unclear, if for any homology class
$e\in H_1(L;\Z)$ there exists an embedding
\begin{gather*}
f:(S^1\times\D,S^1\times\partial\D)
\lra
(\C^2,L),
\end{gather*}
such that the maps $f_{\tau}=f(\tau,.)$ are holomorphic in the usual sense 
for all $\tau\in S^1$ and $[\partial f_{\tau}]=e$.
Such a map (or rather its image) is called an $e${\it -filling} of $L$.

The standard examples of fillings in $\C^2$ are for the Clifford torus
$T_{\pi/2}=\{|z_1|^2+|z_2|^2=1\}$, namely
\begin{gather}
\label{exfill1}
f_1(\tau,z)=
\tfrac{1}{\sqrt{2}}
\big(
z,e^{i\tau}
\big)
\mtqq{and}
f_2(\tau,z)=
\tfrac{1}{\sqrt{2}}
\big(
e^{i\tau},z
\big)
\;.
\end{gather}
Notice that the Maslov index for both cycles is
\begin{gather*}
\mi_{T_{\pi/2}}([\partial f_j])=2
\;,j=1,2.
\end{gather*}

In general it is unknown, if such fillings survive under Hamiltonian
deformations. If one is only interested in the existence of sufficiently 
many holomorphic disks attached to Lagrangian submanifolds such that the 
boundary maps represent certain homology classes, the following relative 
Gromov invariant type notion is useful.

We will say that $L\subset\complex{2}$ admits an $e${\it -quasi-filling} for a homology class 
$e\in H_1(L;\Z)$,
if for all points $p\in L$ there exists a holomorphic disk
\begin{gather*}
u_p:(\D,\partial\D,1)\lra(\C^2,L,p)
\end{gather*}
such that $[\partial u_p]=e$.
{F}or example the following disk families
\begin{gather}
\label{exfill2}
f_a(\tau,z)=
\tfrac{1}{\sqrt{2}}
\left(
z,
e^{i\tau}\frac{z-a}{1-\bar{a}z}
\right)
\mtq{for}
a\in\C\;,\;\abs{a}<1
\end{gather}
which satisfy
\begin{gather*}
\mi_{T_{\pi/2}}([\partial f_a])=4
\end{gather*}
define quasi-fillings for the Clifford torus with respect to the product 
class of the defining circles.
As all disks $f_a(\tau,\D)$ intersect in the singularity $(0,a)$ this is 
not a filling.

We remark that up to parameterizations, the disk families in \eqref{exfill1} 
and \eqref{exfill2}
constitute a complete list of all
unparameterized holomorphic disks attached to the Clifford torus.

It seems that the existence of quasi-fillings depends on the Hamiltonian 
isotopy class as the following statement due to Gromov 
\cite{gro} (\cf
\cite[Proposition 4.1.A]{4knots} and \cite{relgw}) suggests.

\begin{theorem}[{{\bf Gromov}}]\label{twodisks}
Let $L$ be an embedded Lagrangian torus in $\C^2$ Hamiltonian isotopic to the 
Clifford torus $T_{\pi/2}$ via $\varphi$.
Then for $j=1,2$ the monotone torus $L$ admits a $\varphi_*e_j$-quasi-filling,
where $e_1,e_2$ are the generators of $H_1(T_{\pi/2};\Z)$ corresponding to 
the defining circles.
\end{theorem}

Hence, any Lagrangian torus Hamiltonian isotopic to the Clifford torus 
carries two different holomorphic disk families with Maslov index $2$.

This theorem can be used to show that not all embedded Lagrangian tori
are symplectomorphic. We will now give additional examples of monotone
tori. Let $z:S^1\to\complex{*}$ be a closed curve. Using complex-symplectic
coordinates on $\complex{2}$ we obtain a Lagrangian torus through the map
$F:\T^{2}\to\complex{2}$,
$$F(\phi,\psi):=(z(\phi)\cos\psi,z(\phi)\sin\psi).$$
The embedded Lagrangian torus generated by the curve $z(\phi)=e^{i\phi}$
is symplectomorphic to the Clifford torus
$T_{\pi/2}=\set{\abs{z_1}^2+\abs{z_2}^2=1}$ by the $\su_2$-transformation 
matrix
\begin{gather*}
S=
\tfrac{1}{\sqrt{2}}
\left(
\begin{array}{cc}
1&i\\
i&1
\end{array}
\right)
\;.
\end{gather*}
A second example is given by what we call the {\it Chekanov torus} 
$T^2_1(\pi/2)$, see \cite{laglos}.
In this case consider the curve in the open right half plane
\begin{gather*}
z(\phi)=
e^{i\kappa\cos(\phi)}
\cos(\phi)+
i\kappa
e^{-i\kappa\cos(\phi)}
\sin(\phi)
\;,
\mtq{where}
\kappa=\tfrac{1}{\sqrt{2}}
\;.
\end{gather*}
This curve $z$ is embedded with $\{z\}\cap\{-z\}=\emptyset$. The Chekanov torus 
$T^2_1(\pi/2)$ is the embedded equivariant torus generated by $z$, 
see \cite{relgw}.
Equivalently, we mention the {\it Eliashberg-Polterovich torus} 
$L_{\Gamma}^0$ considered in \cite{4knots}, which is symplectically isomorphic 
to the Chekanov torus.
Denote by $\Gamma$ the boundary curve of an embedded disk $D$ in $\C$ such that
$\bar{D}$ does not contain the origin.
By definition $L_{\Gamma}^0$ is the intersection of the sub-level set of
$\Gamma$ under
the function $(z_1,z_2)\mapsto z_1z_2$ with the sub-level set of zero 
under the Hamiltonian
function $H(z_1,z_2)=\abs{z_1}^2-\abs{z_2}^2$.
To see the equivariant picture of this, let $w:S^1\lra\C^*$ be a regular 
parameterisation of $\Gamma$ and choose a branch of square-root defined 
near $D$. Then
\begin{gather}\nonumber
\label{elpolunique}
L_{\Gamma}^0=
\setlr{
\left(
\sqrt{w(\phi)}
e^{it},
\sqrt{w(\phi)}
e^{-it}
\right)
}_{\phi,t\in S^1}
\;.
\end{gather}
Under the $\uni_2$-transformation matrix
\begin{gather}\nonumber
\label{unimat}
U=
\tfrac{1}{\sqrt{2}}
\left(
\begin{array}{cc}
1&i\\
1&-i
\end{array}
\right)=
\left(
\begin{array}{cc}
1&0\\
0&-i
\end{array}
\right)S
\end{gather}
the torus $EQ_{z^{\pt}}$ corresponding to the curve $z^{\pt}=\sqrt{2w}$ 
is mapped to $L_{\Gamma}^0$.
Hence, the torus $EQ_{z^{\pt}}$ is symplectomorphic to $L_{\Gamma}^0$.

The reason why Chekanov and Eliashberg-Polterovich considered these tori is 
that they were the first examples of monotone Lagrangian tori not 
symplectomorphic to the Clifford torus.
Both tori, $T_1^2(\pi/2)$ and $L_{\Gamma}^0$, (which are symplectomorphic) 
are so-called {\sl exotic tori}. The exoticness follows from
Theorem \ref{twodisks}, because $L_{\Gamma}^0$ allows only one holomorphic disk 
family with Maslov index $2$.
{F}urther, for the tori $L_{\Gamma}^0$ we have generating cycles
\begin{gather}
\label{class}
c(\phi)=
\big(
\sqrt{w(\phi)},\sqrt{w(\phi)}
\big)
\end{gather}
for $\phi\in S^1$.
In contrast to Theorem \ref{twodisks} we have the following statement, 
which can be obtained via similar methods (see \cite{relgw}):

\begin{theorem}
Let $L$ be an embedded Lagrangian torus in $\C^2$ symplectomorphic to the 
exotic torus $L_{\Gamma}^0$ via $\varphi$.
Then the monotone torus $L$ admits a $\varphi_*e$-quasi-filling,
where the homology class $e\in H_1(L_{\Gamma}^0;\Z)$ is represented by 
\eqref{class}.
\end{theorem}

An important question is if the existence of $e_i$-quasi-fillings of an
embedded monotone Lagrangian torus $L\subset\complex{n}$ for a full set of 
generators $e_1,\dots, e_n\in H_1(L;\integer{})$ can be used to detect
the type of singularity that it will develop under the mean curvature flow.
We will see that this is not the case (note that the basis of generators is 
in general not canonical).

Essentially, our article consists of two parts. In the first part we 
will derive some results for the Lagrangian mean curvature flow in general, 
whereas in the second part we restrict our attention to an equivariant 
situation to find examples of type-1 and type-2 singularities of monotone
Lagrangian submanifolds $L\subset\complex{n}$ and in particular of embedded, 
monotone Lagrangian tori $\T^{2}\subset\complex{2}$. The equivariant situation
is a generalization of the one discussed above and can be described as 
follows.

Let $L:=S^1\times S^{n-1}$ and assume that

\begin{eqnarray}
z_0:S^1\to\complex{}^*\nonumber\\
z_0(\phi):=u_0(\phi)+iv_0(\phi)\nonumber
\end{eqnarray}

is some closed immersed curve in $\complex{}^*=\complex{}\setminus\{0\}$. 
The equivariant Lagrangian immersions $L_0=F_0(S^1\times S^{n-1})$ 
that we study in the second part of this article, are of the following form. 
Define an immersion $F_0:S^1\times S^{n-1}\to\complex{n}$ by

\begin{eqnarray}
{F}_0(\phi,x):=\bigl(u_0(\phi)G(x),v_0(\phi)G(x)\bigr),\nonumber
\end{eqnarray}

where we assume that the complex structure $J$ is acting on $\complex{n}$ by
$$J(x^1,\dots,x^n,y^1,\dots,y^n):=(-y^1,\dots,-y^n,x^1,\dots,x^n)$$
and that $G:S^{n-1}\to\real{n}$ is the standard embedding of the sphere
of radius $1$ in $\real{n}$. In particular, in
dimension $n=2$ we get nice examples of Lagrangian tori in $\complex{2}$.

Since the mean curvature flow is isotropic, the condition to be equivariant
is preserved under the flow. This implies that in this situation the flow
will be determined by the flow of the corresponding {\sl profile
curves} $\gamma_t:=z_t(S^1)$,  $z_t:=u_t+iv_t:S^1\to\complex{}^*$ induced by 

$$F_t(\phi,x)=\bigl(u_t(\phi)G(x), v_t(\phi)G(x)\bigr).$$

As will be shown in Section \ref{sec flow}, the equivariant Lagrangian mean curvature flow
can be reduced to the following flow of the profile curves in 
$\complex{}^*$

\begin{gather}\label{ecsf}
\dt z=-f\nu,\tag{$*$}\\
z(\cdot,0)=z_0(\cdot),\notag\\
\text{with }f:=k+(n-1)\frac{\langle z,\nu\rangle}{|z|^2}\,.\label{equa 0}
\end{gather}

Here, $\nu$ denotes the outward pointing unit normal along the curve 
$\gamma_t:=z_t(S^1)$, $z_t(\phi)=z(\phi,t)$ and $k$ is the curvature 
of $\gamma_t$.

The reason to consider equivariant solutions is to get simple, non-trivial
examples of monotone solutions that develop type-1 or type-2 singularities.
There are essentially two classes of closed embedded profile curves, those
who contain the origin in their interior and those who do not. In dimension
$n=2$ the Clifford torus belongs to the first class and the exotic tori to 
the second class. In particular, these classes
are not symplectomorphic to each other.

\noindent
{F}or the first class of profile curves we will show 
that neither starshapedness, embeddedness, monotonicity or
a natural convexity assumption are sufficient conditions for the
contraction to a single point (which would have to be the origin in view of 
the point-symmetry). In particular, even under geometrically very restrictive
conditions we will encounter type-2 singularities. For the second class
of curves the first author proves in his PhD thesis \cite{groh}
that they contract in finite time to a single point 
$p\in\complex{*}$, if $f>0$.

On the other hand we can prove fairly general regularity theorems. For example,
one of our main results is:

\begin{theorem}\label{main 1}
Let $F:S^1\times S^{n-1}\times[0,T_{\text{sing}})\to\complex{n}$, $n\ge 2$ be a 
smooth solution of the equivariant Lagrangian mean curvature flow 
with corresponding
profile curves $z:S^1\to\complex{*}$.
If the initial profile curve is starshaped
\footnote{In this paper, a curve $\gamma=z(S^1)$ will be called starshaped, 
if $\langle z,\nu\rangle>0$ everywhere on $\gamma$, no matter if $\gamma$ 
is embedded or immersed.}
w.r.t. the origin, then
$$\liminf_{t\to T_{\text{sing}}}r=0,$$
where $[0,T_{\text{sing}})$ is the maximal time interval on which a smooth solution
exists and $r:=|z|$.
\end{theorem}

This theorem states that for $n\ge 2$ the first singularity of a
starshaped curve must occur at the origin. 

Is is interesting to ask, if there exists a class of monotone Lagrangian
submanifolds that develop type-1 singularities. Let us define

\begin{eqnarray}
\sE:=\{L\subset\complex{n}&:&L\text{ is a monotone, equivariant
Lagrangian} \nonumber\\
&&\text{immersion with $|\overrightarrow H|>0$ and such that
the}\nonumber\\
&&\text{profile curve $z:S^1\to\complex{*}$ is starshaped }\}\nonumber
\end{eqnarray}

Obviously, a Lagrangian $L\subset\sE$ satisfies very
strong geometric conditions.
Therefore, it is surprising that even within this class both types 
of singularities occur; to be precise we will prove the following theorem:

\begin{theorem}\label{main 3}~

\begin{itemize}
\item[(i)]
The class $\sE$ is stable under the mean curvature flow and there exist
embedded $L_1, L_2\subset\sE$ such that $L_1$ develops a singularity of type-1
and $L_2$ a singularity of type-2 at the origin of $\complex{n}$.\\

\item[(ii)]
Let $z:S^1\to\complex{*}$ be the initial profile curve of a Lagrangian
$L\subset\sE$ and let
$l\ge 1+4n\omega_0$ be a positive integer, where $\omega_0=\wind_0(z)$ is 
the winding number of $z$ w.r.t. the origin. In addition, suppose 
that $z$ is rotationally symmetric w.r.t. rotations by 
$\frac{2\pi}{l}$, i.e.
$z(\phi)=z(\phi+\frac{2\pi}{l})$, for all $\phi\in[0,2\pi)$. Then 
the corresponding equivariant Lagrangian immersion $L\subset\sE$
develops a type-1 singularity at the origin and the Lagrangian 
submanifold converges
to the origin as $t\to T_{\text{sing}}=\frac{1}{\varepsilon},$
where $\varepsilon$ is the monotonicity constant of $L$. After rescaling,
the profile curves converge smoothly to a smooth limiting curve 
$\tilde\gamma_\infty$
which is one of the self-similarly contracting solutions classified by 
Anciaux in \cite{anciaux}.
\end{itemize}
\end{theorem}

The organization of our article is as follows. In Section \ref{sec gro} we will
derive general results for the Lagrangian mean curvature flow. There we give
proofs of Theorems and Propositions \ref{lemma mon}, \ref{theo main 1}, 
\ref{prop main 1},  \ref{noself} and \ref{theo main 2}.

In Section \ref{sec equi} we discuss the equivariant situation and prove that
a number of geometric conditions are preserved under the equivariant flow. 
Finally, in the last section we give the proofs for Theorems \ref{main 1} 
and \ref{main 3}.

\section{Mean curvature flow and monotonicity of Lagrangian submanifolds}\label{sec gro}
As in the introduction, let
$F:L\times[0,T_{\text{sing}})\to\complex{n}$ be a maximal smooth solution
of 
\begin{gather}
\dt F=\overrightarrow H\notag\\
F(\cdot,0)=F_0\notag
\end{gather}
for some compact, smooth manifold $L$.
Moreover, for a Lagrangian immersion $F:L\to\complex{n}$ we set
$$F_\alpha:=\frac{\partial F}{\partial y^\alpha},$$
$$F_{\alpha\beta}=\frac{\partial^2 F}{\partial y^\alpha\partial y^\beta}$$ 
and
$$\nu_\alpha:=JF_\alpha,$$
where $\left(y^\alpha\right)_{\alpha=0,\dots,n-1}$ are local coordinates for
$L$. Then, by the Lagrangian condition, the vectors $\nu_\alpha$ 
are normal along $F(L)$.
The induced metric $\dd g\alpha\beta$ and the second fundamental tensor
$\ddd h\alpha\beta\gamma$ on $L$ are given by 
$\dd g\alpha\beta=\langle F_\alpha,F_\beta\rangle$ resp. 
$\ddd h\alpha\beta\gamma=\langle\nu_\alpha, \dd F\beta\gamma\rangle$. 
In local coordinates, the components $H_\alpha$ of the mean curvature
form are given by $H_\alpha=\uu g\beta\gamma\ddd h\alpha\beta\gamma$,
where $\uu g\beta\gamma$ denotes the inverse of $\dd g\beta\gamma$.
The mean curvature vector can be expressed as 
$\overrightarrow H=\uu g\beta\gamma H_\beta\nu_\gamma$.
The following Lemma is well-known (cf. \cite{smoczyk 5}).

\begin{lemma}
Under the mean curvature flow in euclidean space, the mean curvature form
evolves according to the evolution equation

\begin{equation}\label{evol mean}
\dt H=dd^\dagger H,
\end{equation}

where $d^\dagger$ denotes the negative adjoint of $d$ w.r.t. the induced
metric on $L_t$.
\end{lemma}

In addition, as we will now prove, the monotonicity
is preserved, if one allows the monotonicity constant to vary with time.
This was already stated in Proposition \ref{lemma mon}.

\peroof{\ref{lemma mon}}
We need the evolution equation for $F^*\lambda$. One easily observes
$$\left(F^*\lambda\right)_\alpha=\frac{1}{2}\langle JF,F_\alpha\rangle.$$
Then

\begin{eqnarray}
2\dt \left(F^*\lambda\right)_\alpha
&=&\langle J\overrightarrow H,F_\alpha\rangle+
\langle JF,\nabla_\alpha\left(H^\gamma\nu_\gamma\right)\rangle\nonumber\\
&=&-H_\alpha+\nabla_\alpha H^\gamma\langle JF,\nu_\gamma\rangle-H^\gamma h^\beta_{\alpha\gamma}\langle JF,F_\beta\rangle.
\nonumber
\end{eqnarray}

To proceed, we need an expression for $dd^\dagger(F^*\lambda)=d\left(\nabla^\alpha
\left(F^*\lambda\right)_\alpha\right)$.
We compute
$$2\nabla_\alpha(F^*\lambda)_\beta=\langle\nu_\alpha,F_\beta\rangle+h^\gamma_{\alpha\beta}
\langle JF,\nu_\gamma
\rangle=h^\gamma_{\alpha\beta}\langle JF,\nu_\gamma\rangle$$
and
$$2d^\dagger(F^*\lambda)=H^\gamma\langle JF,\nu_\gamma\rangle.$$
Taking the exterior derivative yields
$$2(dd^\dagger(F^*\lambda))_\alpha=\nabla_\alpha H^\gamma\langle JF,\nu_\gamma
\rangle+H_\alpha-
H^\gamma h^\beta_{\alpha\gamma}\langle JF,F_\beta\rangle.$$
Combining this with the above expression for $2\dt(F^*\lambda)_\alpha$ we obtain

\begin{equation}\label{evol lambda}
\dt F^*\lambda=dd^\dagger(F^*\lambda)-H.
\end{equation}

This and the evolution equation (\ref{evol mean}) imply

\begin{eqnarray}
\frac{d}{dt}\left(F^*\lambda-\frac{1}{\varepsilon_0}(1-\varepsilon_0t)H\,\right)
&=&dd^\dagger\left(F^*\lambda-\frac{1}{\varepsilon_0}(1-\varepsilon_0t)H\,\right).
\label{evol mon}
\end{eqnarray}

By the monotonicity at $t=0$ we must have $[F^*\lambda]=\frac{1}{\varepsilon_0}[H]$
at $t=0$. The evolution equation above guarantees that the cohomology
class of $F^*\lambda-\frac{1}{\varepsilon_0}(1-\varepsilon_0t)H$ does not
change. This proves 
$$[F^*\lambda]=\frac{1}{\varepsilon_0}(1-\varepsilon_0t)[H]$$
for all $t\in[0,T_{\text{sing}})$. Now as long as in addition
$t\in[0,\frac{1}{\varepsilon_0})$ we obtain
$$[H]=\frac{\varepsilon_0}{1-\varepsilon_0t}[F^*\lambda]$$
with
$$\frac{\varepsilon_0}{1-\varepsilon_0t}>0$$
which proves the lemma.
\endperoof

We are now ready to prove Theorem \ref{theo main 1}.

\pproof{\ref{theo main 1}}
Suppose $T_{\text{emb}}>\frac{1}{\varepsilon}$. 
Then $T_{\text{sing}}>\frac{1}{\varepsilon}$ and $L_{\frac{1}{\varepsilon}}$
is embedded. From Proposition
\ref{lemma mon} we conclude $[F^*\lambda]=0$ at $t=\frac{1}{\varepsilon}$. 
By Theorem \ref{gromov}, $L_{\frac{1}{\varepsilon}}$ 
cannot be embedded. This contradiction proves
$T_{\text{emb}}\le\frac{1}{\varepsilon}$. 
\endpproof

\peroof{\ref{prop main 1}}
Suppose $F:L\times[0,T_{\text{sing}})\to\complex{n}$ 
is a solution of the Lagrangian mean curvature flow and that
$L_0=F(L,0)$ is monotone with monotonicity constant $\varepsilon$. By Proposition
\ref{lemma mon} we then have
$[F^*\lambda]=\frac{1-\varepsilon t}{\varepsilon}[H]$ for any
$t\in[0,T_{\text{sing}})$.
{F}or the rescaled flow $\widetilde F=\left(2(T_{\text{sing}}-t)\right)^{-\frac{1}{2}}$ we get
$$F^*\lambda=2(T_{\text{sing}}-t)\widetilde F^*\lambda,\quad H=\widetilde H,$$
so that
\begin{equation}\label{tri}
[\widetilde F^*\lambda]
=\frac{1-\varepsilon t}{2\varepsilon(T_{\text{sing}}-t)}[\widetilde H]
\end{equation}
On the other hand, it is well known that a variation of Lagrangian immersions
$\widetilde F_s:L\to\complex{n}$ is Hamiltonian, if and only if 
$$\left[\widetilde F^*\left(\omega\left(\cdot,\frac{d}{ds}\widetilde F\right)\right)\right]=0.$$
{F}or the rescaled flow we have 
$$\left[\widetilde F^*\left(\omega\left(\cdot,\frac{d}{ds}\widetilde F\right)\right)\right]
=\left[\widetilde H-2\widetilde F^*\lambda\right]$$ 
and in view of (\ref{tri}) this vanishes, if and only if $T_{\text{sing}}=\frac{1}{\varepsilon}$.
\endperoof

\pproof{\ref{noself}}
By assumption the mean curvature form satisfies
$$H=d\alpha$$
with 
$$H=\varepsilon F^*\lambda.$$
This implies
$$d^\dagger H=\Delta\alpha=\varepsilon d^\dagger F^*\lambda=
\varepsilon\langle d\alpha,ds\rangle$$
with $s:=\frac{|F|^2}{4}$. We multiply this equation 
with $\alpha e^{-\varepsilon s}$ and integrate on $L$. 
Note, that the integral exists in view of the curvature estimate 
for $|\nabla^k\alpha|^2$. This gives
$$\int_L\alpha \Delta\alpha e^{-\varepsilon s}d\mu=\int_L\varepsilon\alpha
\langle d\alpha,ds\rangle e^{-\varepsilon s}d\mu.$$
By assumption $|\nabla^k\alpha|^2\le ce^{\sigma s}$ with 
$\sigma <\varepsilon$
so that we may integrate by parts in the first integral without getting
boundary terms at infinity. This yields
$$\int_L|H|^2e^{-\varepsilon s}d\mu=0,$$
hence $H=\varepsilon F^*\lambda=0$ which means that $L$ is part of a
minimal cone.
\endpproof

\pproof{\ref{theo main 2}}

Suppose $L_t$ develops an embedded type-1 singularity at the origin and that
$\tilde L_\infty$ is an embedded limiting Lagrangian of the rescaled flow 
$\tilde L_s, s\in[s_0,\infty)$ with $s_0:=-\frac{1}{2}\log T_{\text{sing}}$. 
Then by Huisken's
result we know $\tilde H_\infty=2\tilde\lambda_\infty$. 

If $[\widetilde H_\infty]=0$, then in view of Theorem \ref{noself} 
$\tilde L_\infty$ would be part of a minimal cone. On the other hand
it is well-known that a type-1 blow-up of a compact submanifold
is smooth, complete and not totally geodesic
but any smooth, complete part of a minimal cone is totally geodesic.
So we must have $[\widetilde H_\infty]\neq 0$. 

$[\widetilde H_\infty]\neq 0$ and 
$\widetilde H_\infty=2\widetilde\lambda_\infty$ implies that
there exists a closed embedded curve 
$\widetilde\gamma\subset\widetilde L_{\infty}$ 
such that 
$\int_{\widetilde\gamma}\tilde\lambda_\infty\neq 0$. Choose a ball $B(0,R)\subset\complex{n}$ 
with $\widetilde\gamma\subset B(0,R)$ for some $R>0$. Since $\widetilde L_\infty$ is
an embedded limiting Lagrangian submanifold, there exist homology classes
$\widetilde e_\infty\in H_1(\widetilde L_\infty;\integer{})$, 
$e\subset H_1(L;\integer{})$,
a closed embedded curve $\widetilde\gamma_\infty\in\widetilde e_\infty$ and a 
sequence of closed embedded curves 
$\gamma_t\in e$ such that $\widetilde\gamma_\infty\subset B(0,R)$, 
$\int_{\widetilde\gamma_\infty}\tilde\lambda_\infty\neq 0$,
$\widetilde\gamma_s:=\widetilde F(\gamma_t,s(t))\subset B(0,R)$. 
Note, that the type-1 assumption implies that $|\nabla^k\widetilde A|^2$ is uniformly
bounded for $k\ge 0$ and all $s$. 
So, using the compactness of $B(0,R)$, we may assume that $\widetilde\gamma_\infty$
is the limiting curve of $\widetilde\gamma_s$ and that 
$$\lim_{s\to\infty}\int_{\widetilde\gamma_s}\tilde\lambda_s=\int_{\widetilde\gamma_\infty}\tilde\lambda_\infty,$$
$$\lim_{s\to\infty}\int_{\widetilde\gamma_s}\widetilde H_s=\int_{\widetilde\gamma_\infty}\widetilde H_\infty.$$
A direct computation shows that under the
rescaled flow the $1$-form $\widetilde H_s-2\widetilde\lambda_s$ satisfies the
evolution equation
$$\frac{d}{ds}\left(\widetilde H_s-2\widetilde\lambda_s\right)
=dd^\dagger\left(\widetilde H_s-2\widetilde\lambda_s\right)
+2\left(\widetilde H_s-2\widetilde\lambda_s\right).$$
Since $\widetilde \gamma_s=\widetilde F(\gamma_t,s(t))$ and 
$\gamma_t\in e,\forall\,t$
we conclude
$$\int_{\widetilde\gamma_s}\left(\widetilde H_s-2\widetilde\lambda_s\right)
=T_{\text{sing}}e^{2s}\int_{\widetilde\gamma_{s_0}}\left(\widetilde H_{s_0}-2\widetilde\lambda_{s_0}\right).$$
On the other hand
$$\lim_{s\to\infty}\int_{\widetilde\gamma_s}\left(\widetilde H_s-2\widetilde\lambda_s\right)=0$$
implies
$$\int_{\widetilde\gamma_{s_0}}\left(\widetilde H_{s_0}-2\widetilde\lambda_{s_0}\right)=0.$$
Since $L_0$ is monotone we can use equation (\ref{tri}) and $[\widetilde H_s]=[\widetilde H_{s_0}]$
to obtain
$$\left(1-\frac{1}{\varepsilon T_{\text{sing}}}\right)
\int_{\widetilde\gamma_{s_0}}\widetilde H_{s_0}=0.$$
So either $T_{\text{sing}}=\frac{1}{\varepsilon}$ or
$\int_{\widetilde\gamma_{s_0}}\widetilde H_{s_0}=0.$
Since $\frac{d}{ds}\int_{\widetilde\gamma_s}\widetilde H_s=0$ and 
$\lim_{s\to\infty}\int_{\widetilde\gamma_s}\widetilde H_s
=\int_{\widetilde\gamma_\infty}\widetilde H_\infty$ we obtain
$\int_{\widetilde\gamma_\infty}\widetilde H_\infty=0$ in the last case. 
This is a contradiction
to $\widetilde H_\infty =2\widetilde\lambda_{\infty}$ and
$\int_{\widetilde\gamma_\infty}\widetilde\lambda_\infty\neq 0$. So 
$T_{\text{sing}}=\frac{1}{\varepsilon}$ and we are done.
\endpproof

\section{Mean curvature flow of equivariant Lagrangian immersions of 
$S^1\times S^{n-1}$}\label{sec equi}
\subsection{The equivariant flow}\label{sec flow}~

We want to show that the Lagrangian mean curvature flow of equivariant Lagrangian
immersions takes the form given in (\ref{ecsf}). Therefore, assume
$$F(\phi,x)=\bigl(u(\phi)G(x),v(\phi)G(x)\bigr)$$
is an equivariant Lagrangian immersion as above. We will denote the coordinate
on $S^1$ by $\phi$ and local coordinates on $S^{n-1}$ by $x^1,\dots,x^{n-1}$.
Latin indices $i, j, k, ...$ will be in the range between $1$ and $n-1$, 
whereas greek indices
$\alpha,\beta, ...$ are taken between $0$ and $n-1$. In particular, we define
coordinates $y^\alpha$ on $S^1\times S^{n-1}$ by $y^0:=\phi, y^i:=x^i$
for all $i\in\{1,\dots,n-1\}$. Doubled indices will be summed according to
the Einstein convention, i.e. latin indices from $1$ to $n-1$ and greek
indices from $0$ to $n-1$.

We want to compute the
metric and second fundamental form, in particular the mean curvature of these
equivariant Lagrangian immersions. To this end, let us denote any partial derivative
of $u, v$ w.r.t. $\phi$ by a prime and in addition we set
$G_i:=\frac{\partial G}{\partial x^i}$ and
$G_{ij}=\frac{\partial^2 G}{\partial x^i\partial x^j}$. 
With this notation we get 
\begin{gather}
{F}_0=(u' G,v' G)\,, F_i=(uG_i,vG_i),\nonumber\\
\nu_0=(-v' G,u' G)\,, \nu_i=(-vG_i,uG_i),\nonumber\\
{F}_{00}=(u''G,v''G)\,,F_{0 i}=(u' G_i,v' G_i)\,,F_{ij}=(uG_{ij},vG_{ij}).
\nonumber
\end{gather}
The induced metric $\dd g\alpha\beta$ and the second fundamental tensor
$\ddd h\alpha\beta\gamma$ on $L$ are given by 
$\dd g\alpha\beta=\langle F_\alpha,F_\beta\rangle$ resp. 
$\ddd h\alpha\beta\gamma=\langle\nu_\alpha, \dd F\beta\gamma\rangle$. 
The standard metric on
$S^{n-1}$ will be denoted by $\dd\sigma ij$.
Thus
\begin{gather}
\dd g00=(u')^2+(v')^2\,,\,\dd g0i=0\,,\,\dd gij=(u^2+v^2)\dd\sigma ij\nonumber
\end{gather}
and
\begin{gather}
\ddd h000=u'v''-v'u''\,,\,
\ddd h00i=0\,,\,
\ddd h0ij=(uv'-vu')\dd\sigma ij\,,\,
\ddd hijk=0\,.\nonumber
\end{gather}
{F}or the mean curvature form $H_\alpha=\uu g\beta\gamma\ddd h\alpha\beta\gamma$
we obtain
\begin{eqnarray}
\d H0&=&(n-1)\frac{uv'-vu'}{u^2+v^2}+\frac{u'v''-v'u''}{(u')^2+(v')^2},
\label{mcf a}\\
\d Hi&=&0.\label{mcf b}
\end{eqnarray}
In particular the mean curvature vector $\overrightarrow H$ is given by
\begin{eqnarray}
\overrightarrow H&=&\uu g\alpha\beta H_\alpha\nu_\beta\nonumber\\
&=&\uu g00H_0\nu_0\nonumber\\
&=&\frac{1}{(u')^2+(v')^2}
\left((n-1)\frac{uv'-vu'}{u^2+v^2}
+\frac{u'v''-v'u''}{(u')^2+(v')^2}\right)(-v'G,u'G).\nonumber
\end{eqnarray}
We want to rewrite the last equation in terms of objects on the curve
$\gamma\subset\complex{}^*$. If we orient $\gamma$ in the usual way, then
the outer unit normal $\nu$ along $\gamma$ is given by
$$\nu=-J\left(\frac{z'}{|z'|}\right).$$
This implies
$$\nu=\frac{1}{\sqrt{(u')^2+(v')^2}}\binom{v'}{-u'}$$
and
$$\langle z,\nu\rangle=\frac{uv'-vu'}{\sqrt{(u')^2+(v')^2}}.$$
The curvature $k$ of the curve $\gamma$ is determined by
$$k=-\frac{1}{|z'|^2}\langle z'',\nu\rangle=\frac{u'v''-v'u''}
{\bigl((u')^2+(v')^2\bigr)^{3/2}}.$$
Then these equations imply
\begin{equation}\label{mean}
\overrightarrow H=-\left(k+(n-1)\frac{\langle z,\nu\rangle}{|z|^2}\right)
\frac{1}{\sqrt{(u')^2+(v')^2}}(v'G,-u'G).
\end{equation}

If we project the Lagrangian mean curvature flow (MCF) to
the ``north-pole plane'' $((1,0,\dots,0),(1,0,\dots,0))\subset\complex{n}$,
then the flow 
$$\dt F=\overrightarrow H=(\dt u\ G,\dt v\ G)$$ 
induces the flow
$$\dt z=-\left(k+(n-1)\frac{\langle z,\nu\rangle}{|z|^2}\right)\nu$$
which is equation (\ref{ecsf}). 

\subsection{Singularities}~

In case of an equivariant
Lagrangian immersion as above, one easily observes that the 
second fundamental tensor $A$ on $L$ can be completely expressed in terms
of the curvature $k$ of the curves and 
$\frac{\langle \nu,e_r\rangle}{r}$
so that we obtain 
\begin{lemma}\label{lemma long}
Let $\gamma\subset\complex{*}$ be a closed curve that evolves under
the flow (\ref{ecsf}).
As $t$ tends to $T_{\text{sing}}$ we must either have
$$\limsup_{t\to T_{\text{sing}}}|k|=\infty\quad\text{or}\quad\liminf_{t\to T_{\text{sing}}}r=0.$$
Moreover, if $n>1$ and if there exists $\tilde T\in(0,T_{\text{sing}}]$ such that
$$\liminf_{t\to \tilde T}r=0,$$
then $\tilde T=T_{\text{sing}}$.
\end{lemma}
\begin{proof}
At a point $\phi_0\in S^1$, where $r(\phi_0)=\min_{\phi\in S^1}r(\phi)$
we have $z=\pm r\nu$ 
and it follows that at such a point $\frac{|\langle\nu,e_r\rangle|}{r}=
\frac{1}{r}$. If $n>1$ and
$r$ tends to zero somewhere, then the corresponding equivariant Lagrangian 
map fails to be an immersion as $t\to\tilde T$.
\end{proof}
In view of the last lemma we distinguish the following three cases:
\begin{gather}
\limsup_{t\to T_{\text{sing}}}|k|=\infty\quad\text{ and }
\quad\liminf_{t\to T_{\text{sing}}}r=0,\tag{C1}\\
\limsup_{t\to T_{\text{sing}}}|k|=\infty\quad\text{ and }
\quad\liminf_{t\to T_{\text{sing}}}r>0,\tag{C2}\\
\limsup_{t\to T_{\text{sing}}}|k|<\infty\quad\text{ and }
\quad\liminf_{t\to T_{\text{sing}}}r=0.\tag{C3}
\end{gather}
\begin{remark}
A priori we cannot exclude a singularity of type $(C3)$ $($except in the
case $n=1)$. One can prove that in case $(C3)$
the curves converge in the $C^\infty$-topology to a smooth limit curve 
$\gamma$ passing through the origin as $t\to T_{\text{sing}}$. The singularity for the
corresponding Lagrangian immersions then occurs, because in the limit
$t\to T_{\text{sing}}$ the equivariant maps $F:S^1\times S^{n-1}\to\complex{n}$ fail to
be an immersion. We thus call a singularity of type $(C3)$ a fake
singularity. The next lemma states that the quantity 
$\frac{\langle\nu,e_r\rangle}{r}$ and therefore
also the equivariant Lagrangian mean curvature flow would still make 
sense for those curves.
\end{remark}
\begin{lemma}\label{lemma 1}
Let $z:S^1\to\complex{}$ be a smooth regular curve. Assume $\phi_0\in S^1$
is a point with $z(\phi_0)=0$. Then we have
$$\lim_{\phi\to\phi_0}\frac{\langle z(\phi),\nu(\phi)\rangle}{|z(\phi)|^2}
=\frac{1}{2}k(\phi_0),$$
where $k(\phi_0)$ denotes the curvature of the curve at $z(\phi_0)=0$ and
$\nu(\phi)$ is the inward pointing unit normal.
\end{lemma}
\begin{proof}
Let $(x,y)$ be cartesian coordinates for $\complex{}$ such that at the origin
the $x$-axis is tangent to our curve $\gamma:=z(S^1)$. Then locally, say
on an interval $(-\epsilon,\epsilon)$, we can express $\gamma$ as a graph
over the $x$-axis, i.e. w.l.o.g. we may assume that after a reparameterization
our curve is locally given by $z(x)=(x,y(x))$ with a smooth function 
$y:(-\epsilon,\epsilon)\to\real{}$. It is sufficient to prove the lemma for
real analytic functions $y$. Let
$$y(x)=\sum_{k=0}^\infty\frac{x^k}{k!}y^{(k)}(0)$$
be the Taylor expansion of $y$ around $0$.
{F}or $x\neq 0$ we have
$$\frac{\langle z,\nu\rangle}{|z|^2}=\frac{xy'-y}{(x^2+y^2)\sqrt{1+(y')^2}}
=\frac{\frac{y'}{x}-\frac{y}{x^2}}{\left(1+\left(\frac{y}{x}\right)^2\right)
\sqrt{1+(y')^2}}.$$ 
Since $y(0)=0$ we see that 
$$\lim_{x\to0}\frac{y^2}{x^2}=(y'(0))^2.$$
If we take into account the Taylor expansion we obtain
$$\frac{1}{x^2}(y-xy')=\sum_{k=2}^\infty x^{k-2}y^{(k)}(0)\left(\frac{1}{k!}
-\frac{1}{(k-1)!}\right)$$
and this implies
$$\lim_{x\to 0}\frac{1}{x^2}(y-xy')=-\frac{1}{2}y''(0)$$
so that
$$\lim_{x\to 0}\frac{\langle z(x),\nu(x)\rangle}{|z(x)|^2}=\frac{y''(0)}
{2\left(1+(y'(0))^2\right)^{3/2}}=\frac{1}{2}k(0).$$
This proves the lemma.
\end{proof}

\begin{remark}
One has to be a bit careful here. If $z:S^1\times[0,T_{\text{sing}})\to\complex{}^*$ is a 
smooth family of curves with $\liminf_{t\to T_{\text{sing}}}r=0$, then clearly 
$$\limsup_{t\to T_{\text{sing}}}\frac{|\langle\nu,e_r\rangle|}{r}=\infty.$$ 
But if a smooth limit curve (passing through
the origin) exists, then lemma \ref{lemma 1} does not prove that the
curvature $k$ blows up at the origin (but $|A|$ does for 
$t\to T_{\text{sing}}$). 
This shows that (\ref{ecsf}) would also be well defined for a smooth curve 
passing through the origin and a solution of (\ref{ecsf}) might exist longer
than the corresponding solution for (MCF). On the other hand, it is likely
that fake singularities will not occur at all under the mean curvature flow. 
\end{remark}

\subsection{Embeddedness}~

The next question we address is, if the embeddedness of the profile curve
$\gamma\subset\complex{*}$ will also imply the embeddedness of the
equivariant Lagrangian submanifold associated to $\gamma$.

\begin{lemma}
Let $\gamma\subset\complex{*}$ be an embedded curve given by an embedding
$z:S\sp 1\to\complex{*}$. The image $L:=F(S\sp 1\times S\sp{n-1})$, 
$F(\phi,p)=\Bigl(u(\phi)p,v(\phi)p\Bigr)$ is an embedded Lagrangian 
submanifold, if and only if one of the following conditions is valid: 
(i) the curve $\gamma$ satisfies $\gamma\cap(-\gamma)=\emptyset$, 
(ii) $\gamma$ is point-symmetric, i.e. it satisfies $\gamma=-\gamma$. 
In the  latter case, $F$ is a double covering of $L$.
\end{lemma}
\begin{proof}
Let $(\phi_1,p_1), (\phi_2,p_2)\in S\sp1\times S\sp{n-1}$ to points with
$F(\phi_1,p_1)=F(\phi_2,p_2)$. Then either $p_1=p_2, z(\phi_1)=z(\phi_2)$
which, in view of the embeddedness of $\gamma$ implies $\phi_1=\phi_2$,
or $p_1=-p_2,z(\phi_1)=-z(\phi_2)$. So if $\gamma\cap(-\gamma)\neq\emptyset$
and $\gamma\cap(-\gamma)\neq\gamma$, then $L$ cannot be embedded. The rest
of the lemma is clear.
\end{proof}

\begin{lemma}\label{lemma emb}
If the initial curve is embedded, then this remains true for all 
$t\in[0,T_{\text{sing}})$.
\end{lemma}
\begin{proof}
This is well-known for $n=1$. Assume $n>1$
and that the curve touches itself for the first time at $t_0<T_{\text{sing}}$. Let
$\phi_1,\phi_2\in S^1$ be two distinct points with $z(\phi_1)=z(\phi_2)$.
Since $t_0<T_{\text{sing}}$ we must have $|z|>0$ on $S^1$ (Lemma \ref{lemma long}). 
It suffices to prove that
$\langle f(\phi_1)\nu(\phi_1)-f(\phi_2)\nu(\phi2),\nu(\phi_1)-\nu(\phi_2)
\rangle\le0$ 
because this contradicts the strong maximum principle. The curvatures at
$\phi_1$ and $\phi_2$ must satisfy $k(\phi_1)+k(\phi_2)\le 0$. Moreover,
we have $\nu(\phi_1)=-\nu(\phi_2)$, $z(\phi_1)=z(\phi_2)$ and
$\frac{\nur}{r}(\phi_1)=-\frac{\nur}{r}(\phi_2)$. Therefore one has
\begin{eqnarray}
&&\langle f(\phi_1)\nu(\phi_1)-f(\phi_2)\nu(\phi2),\nu(\phi_1)-\nu(\phi_2)
\rangle\nonumber\\
&&=2\bigl(k(\phi_1)+(n-1)\frac{\nur}{r}(\phi_1)
+k(\phi_2)+(n-1)\frac{\nur}{r}(\phi_2)\bigr)\nonumber\\
&&=2\bigl(k(\phi_1)+k(\phi_2)\bigr)\le 0.\nonumber
\end{eqnarray}
\end{proof}

\begin{remark}
The lemma does not imply that the curves stay embedded as $t\to T_{\text{sing}}$. At
a singularity this property might fail. 
\end{remark}

\subsection{Tamed curves}
\begin{definition}
A curve $\gamma\subset\complex{*}$ will be called {\sl tamed}, if the
quantity
$$f:=k+(n-1)\frac{\langle z,\nu\rangle}{|z|^2}$$
is positive for all $p\in\gamma$.
\end{definition}

In the sequel we will set
$$r:=|z|\,,\quad e_r:=\frac{z}{|z|}$$
so that
$$f=k+(n-1)\frac{\langle\nu,e_r\rangle}{r}\,.$$

\begin{lemma}\label{lemma 2}
Let $\gamma\subset\complex{*}$ be a smooth, closed curve. Then there exists 
at least one point $p\in\gamma$ where 
$$f\ge\frac{n}{r}>0.$$
\end{lemma}
\begin{proof}
Let $z:S\sp 1\to\complex{*}$ be a smooth immersion of $\gamma$ and 
$\phi\in S\sp 1$ a point where $r$ assumes its maximum. By definition
of the outer unit normal we have at this point $z=r\nu$. 
Since at $\phi$ we must have $\Delta r\sp 2=\Delta|z|\sp 2\le 0$, the
identity $\lap z=-k\nu$ gives
$$0\ge\lap r^2=2(1-k\langle z,\nu\rangle)=2(1-kr).$$
Therefore 
$$f=k+(n-1)\frac{\langle \nu,e_r\rangle}{r}
\ge\frac{1}{r}+(n-1)\frac{1}{r}
=\frac{n}{r}\,.$$
\end{proof}

{F}rom this and equation (\ref{mean}) follows:
\begin{corollary}\label{coral 1}
A curve $\gamma\subset\complex{*}$ is tamed if and only if the
corresponding equivariant Lagrangian immersion of $S\sp 1\times S\sp{n-1}$
into $\complex{n}$ has a nowhere vanishing mean curvature vector field
$\overrightarrow H$.
\end{corollary}

Next we will compute several evolution equations related to 
(\ref{ecsf}). From \cite{smoczyk 1} we first obtain the equations:
\begin{lemma}\label{lemma 3}
\begin{eqnarray}
\dt d\mu&=&-kfd\mu\label{evol 1}\\
\dt \nu&=&\nabla f\label{evol 2}\\
\dt k&=&\lap f+fk^2\label{evol 3},
\end{eqnarray}
where $d\mu$ denotes the induced volume form on $\gamma$.
\end{lemma}

\begin{lemma}\label{lemma 4}
The evolution equation for $f$ is given by
\begin{equation}\label{evol 4}
\dt f=\lap f+\frac{n-1}{r}\langle\nabla f,e_r\rangle
+f\left(k^2+\frac{n-1}{r^2}\bigl(2\langle\nu,e_r\rangle^2-1\bigr)\right).
\end{equation}
\end{lemma}
\begin{proof}
We have
$$\dt\frac{\langle z,\nu\rangle}{r^2}=\frac{-f+\langle z,\nabla f\rangle}{r^2}
+\frac{2\langle z,\nu\rangle^2f}{r^4}$$
and then with (\ref{evol 3})
$$\dt f=\Delta f+fk^2+(n-1)\left(\frac{-f+\langle z,\nabla f\rangle}{r^2}
+\frac{2\langle z,\nu\rangle^2f}{r^4}\right).$$
\end{proof}

This and the maximum principle implies:
\begin{corollary}\label{coral 2}
If the initial curve $\gamma_0\subset\complex{*}$
is tamed, then it stays tamed for $t\in[0,T_{\text{sing}})$.
\end{corollary}

If we take into account Corollary \ref{coral 1}, then also
\begin{corollary}\label{coral 2b}
If the initial equivariant Lagrangian immersion $L$
satisfies $\overrightarrow H\neq 0$ for all $x\in L$, 
then this remains true for all $t\in[0,T_{\text{sing}})$.
\end{corollary}

\subsection{Monotonicity}~

We want to find the condition when an equivariant Lagrangian immersion
with a profile curve $\gamma\subset\complex{*}$ becomes monotone. Since
we have already computed the mean curvature form in (\ref{mcf a}) and
(\ref{mcf b}), we need to compute $F^*\lambda$ only.

Since
$$\left(F^*\lambda\right)_\alpha=\frac{1}{2}\langle JF,F_\alpha\rangle$$
and 
\begin{eqnarray}
JF&=&(-vG,uG),\nonumber\\
{F}_\alpha&=&(u_\alpha G+uG_\alpha,v_\alpha G+vG_\alpha),\nonumber\\
|G|^2&=&1,\nonumber\\
\langle G,G_\alpha\rangle&=&0\nonumber
\end{eqnarray}

we obtain:
\begin{equation}\nonumber
\left(F^*\lambda\right)_\alpha=\frac{1}{2}
\left(uv_\alpha-vu_\alpha\right),
\end{equation}

i.e.
\begin{eqnarray}
\left(F^*\lambda\right)_0=\frac{1}{2}(uv'-vu'),\label{mon a}\\
\left(F^*\lambda\right)_i=0,\quad\forall\, i=1,\dots, n.\label{mon b}
\end{eqnarray}

{F}rom equations (\ref{mcf a}), (\ref{mcf b}), (\ref{mon a}) and (\ref{mon b})
we see that $[H]=\varepsilon[F^*\lambda]$, if and only if
$$\int_\gamma H=\varepsilon\int_\gamma z^*\lambda.$$
Now, since by (\ref{mcf a}) and (\ref{mon a})
$$H_0=fd\mu,\quad (F^*\lambda)_0=\frac{1}{2}\langle z,\nu\rangle d\mu$$
we obtain

\begin{lemma}\label{monoton equi}
An equivariant Lagrangian immersion is monotone with monotonicity constant
$\varepsilon$, if and only if on the profile curve $z:S^1\to\complex{*}$
we have
\begin{equation}\nonumber
\int_{S^1} fd\mu
=\frac{\varepsilon}{2}\int_{S^1}\langle z,\nu\rangle d\mu.
\end{equation}
\end{lemma}

If we define the symplectic area $A(z)$ of a regular, smooth
curve $z:S^1\to\complex{*}$ by 
$$A(z):=\int_{S^1}z^*\lambda,$$
then the computations above show
$$A(z)=\frac{1}{2}\int_{S^1}\langle z,\nu\rangle d\mu.$$
Moreover, the next lemma gives a geometric interpretation of
$\int_{S^1} fd\mu$ as well.

\begin{lemma}\label{lemma wind}
Suppose $z:S^1\to\complex{*}$ is a regular, smooth curve. Let
$\wind_0(z)$ and $\rot(z)$ denote the winding number of $z$ 
around the origin resp. the rotation number. Then 
\begin{equation}
\int_{S^1} fd\mu=2\pi\left(\rot(z)+(n-1)\wind_0(z)\right).
\end{equation}
\end{lemma}
\begin{proof}
The vector field $Z(z):=z$ satisfies
$$\text{div}\,(Z)
=2\quad\text{and}\quad\text{div}\,\left(\frac{Z}{|z|^2}\right)
=0,\text{ for } z\neq 0.$$
Stokes theorem gives
$$A(z)=\frac{1}{2}\int_{S^1}\langle z,\nu\rangle d\mu$$
and
$$2\pi\wind_0(z)=\int_{S^1}\left\langle\frac{z}{|z|^2},\nu\right\rangle d\mu.$$
Since 
$$f=k+(n-1)\left\langle\frac{z}{|z|^2},\nu\right\rangle$$ 
and
$$\int_{S^1}kd\mu=2\pi\rot(z)$$ 
the proof follows.
\end{proof}

{F}rom Lemma \ref{monoton equi} and Lemma \ref{lemma wind} we immediately 
obtain:
\begin{corollary}
Suppose $z:S^1\to\complex{*}$ is a regular, smooth curve with $A(z)\neq 0$
and $(\rot(z)+(n-1)\wind_0(z))/A(z)>0$. Then the corresponding
equivariant Lagrangian immersion is monotone with monotonicity constant
$$\varepsilon =\frac{2\pi\left(\rot(z)+(n-1)\wind_0(z)\right)}{A(z)}>0.$$
\end{corollary}

\subsection{The symplectic area formula}~

\begin{theorem}\label{theorem 1}
Suppose $z:S^1\times[0,T_{\text{sing}})\to\complex{}^*$ is a smooth solution of (\ref{ecsf}). 
Then the symplectic area $A(z_t)$ of $z_t:=z(\cdot,t):S^1\to\complex{*}$,
satisfies the following equation

\begin{equation}\label{evol 7}
A(z_t)=A(z_0)-2\pi\bigl(\rot(z_0)
+(n-1)\wind_0(z_0)\bigr)t.
\end{equation}

\end{theorem}
\begin{proof}
Since $\dt\langle z,\nu\rangle=-f+\langle\nabla f,z\rangle$ and
$\dt d\mu=-fk d\mu$ we obtain
$$\dt A(z)=\frac{1}{2}\int_{S^1}\bigl(\langle\nabla f,z\rangle-f
-\langle z,\nu\rangle fk\bigr) d\mu$$
and with partial integration
$$\dt A(z)=\frac{1}{2}\int_{S^1}\bigl(-f(1-k\langle z,\nu\rangle)
-f-\langle z,\nu\rangle fk\bigr) d\mu=-\int_{S^1}fd\mu.$$
By Lemma \ref{lemma wind} we thus conclude
$$\dt A(z_t)=-2\pi\left(\rot(z_t)+(n-1)\wind_0(z_t)\right)$$
The rotation number is always invariant under homotopies so that 
$\rot(z_t)=\rot(z_0)$. 
{F}or $n>1$, Lemma \ref{lemma long} implies $|z|>0$ on $[0,T_{\text{sing}})$. 
Thus, in case $n>1$, $\wind_0(z_t)=\wind_0(z_0)$. This proves the theorem.
\end{proof}

\begin{corollary}
Suppose $z:S^1\times[0,T_{\text{sing}})\to\complex{}^*$ is a smooth 
solution of (\ref{ecsf}) and that $z_0:S^1\to\complex{*}$ is an embedding.
Then 
$$T_{\text{sing}}\le\frac{A(z_0)}{2\pi(\rot(z_0)+(n-1)\wind_0(z_0))}.$$ 
\end{corollary}
\begin{proof}
{F}irst note, that this does not follow from Theorem \ref{theo main 1}
since the equivariant Lagrangian submanifold need not be embedded. 
On the other hand
an embedded closed curve satisfies $\rot(z_0)=1$ and $\wind_0(z_0)\in\{0,1\}$
depending on whether the curve encloses the origin or not. Since $A(z_0)$ is
the area of the enclosed region of $z_0$ we see that $A(z_t)>0$ for
all $t\in[0,T_{\text{sing}})$ because the curves stay embedded. The
symplectic area formula above shows that this can be true only as long as
$$T_{\text{sing}}\le\frac{A(z_0)}{2\pi(\rot(z_0)+(n-1)\wind_0(z_0))}.$$ 
Note also, that the equivariant Lagrangian submanifold
for $z_0$ is monotone with monotonicity constant
$$\varepsilon=\frac{2\pi(\rot(z_0)+(n-1)\wind_0(z_0))}{A(z_0)}$$ 
so that we have
$$T_{\text{emb}}=T_{\text{sing}}\le\frac{1}{\varepsilon}.$$
\end{proof}

\subsection{Selfsimilar solutions}~

Like for the monotonicity, the equation for self-similar solutions 
(\ref{monotone}) takes a much simpler form in the equivariant case, namely
\begin{equation}\label{equa 1}
f=\frac{\varepsilon}{2}\langle z,\nu\rangle.
\end{equation}

Henri Anciaux \cite{anciaux} classified all solutions of (\ref{equa 1}).
He noticed that - in contrast to the curve-shortening flow - not all solutions
of (\ref{equa 1}) are convex. On the other hand, in view of Corollary 
\ref{coral 1}, tameness is a much more natural condition than convexity
and by (\ref{equa 1}) any tamed solution is also 
starshaped.
Indeed, as the next theorem shows, all
solutions of (\ref{equa 1}) with $\varepsilon\neq 0$ are tamed and starshaped.

\begin{theorem}
Assume $z$ is a self-similar expanding, or contracting solution of
$$\frac{d}{dt}z=-f\nu,$$
i.e. there exists a constant $\varepsilon\neq 0$ such that
\begin{equation}\nonumber
f=\frac{\varepsilon}{2}\langle z,\nu\rangle.
\end{equation}
Then the quantity
$$p:=fe^{-\frac{\varepsilon}{4}|z|^2}|z|^{n-1}$$
is constant along the curve and the curve is tamed and starshaped.
\end{theorem}
\begin{proof}
The Weingarten equation implies
\begin{equation}\label{weingarten}
\nabla\langle z,\nu\rangle=\frac{1}{2}\,k\nabla|z|^2.
\end{equation}
In view of (\ref{equa 0}) and (\ref{equa 1}) we obtain
\begin{equation}\label{equa 2}
k=\left(1-2\frac{n-1}{\varepsilon|z|^2}\right)f.
\end{equation}
If we take the gradient of both sides in (\ref{equa 1}), then we obtain
\begin{eqnarray}
\nabla f
&=&\frac{\varepsilon}{2}\nabla\langle z,\nu\rangle\nonumber\\
&=&\left(\frac{\varepsilon}{4}-\frac{n-1}{2|z|^2}\right)f\nabla|z|^2,\label{equa 3}
\end{eqnarray}
where we have used (\ref{weingarten}) and (\ref{equa 2}) in the
last step.
We claim that the quantity
$$p:=fe^{-\frac{\varepsilon}{4}|z|^2}|z|^{n-1}$$
is constant on any solution of (\ref{equa 0}). The quantity $p$ is
differentiable whenever $z\neq 0$.
By (\ref{equa 3}) the
gradient of $p$ at those points is
\begin{eqnarray}
\nabla p
&=&\left(\frac{\varepsilon}{4}-\frac{n-1}{2|z|^2}\right)p\nabla|z|^2
-\frac{\varepsilon}{4}\,p\nabla|z|^2+\frac{n-1}{|z|}p\nabla|z|\nonumber\\
&=&0\nonumber.
\end{eqnarray}
The constant $p$ can only be zero, if $f\equiv0$ which is equivalent to
$z$ being a segment of a line passing through the origin. In particular, any 
solution of (\ref{equa 1}) passing through the origin must be a line segment.
\\~\\
{F}or a curve different from a line there exists always one point, where 
$\langle z,\nu\rangle>0$. Consequently, by (\ref{equa 1}), $p$ must be a 
positive constant and it also follows that
any solution  $z:S^1\to\mathbb{C}^*$  of (\ref{equa 1}) is starshaped 
w.r.t. the origin and satisfies $f>0$. 
\end{proof}

\begin{remark}
In particular, by Corollary \ref{coral 1}, the corresponding
self-similarly evolving equivariant Lagrangian immersions in
$\mathbb{C}^n$ all have nowhere vanishing mean curvature vector 
$\overrightarrow H$. Moreover, as has been noted earlier by Anciaux
\cite{anciaux}, for $n=1$ the curves are all convex,
whereas for $n>1$ we have
\begin{eqnarray}
k>0,&& \text{for }|z|>\sqrt{\frac{2(n-1)}{\varepsilon}}\nonumber\\
k=0,&& \text{for }|z|=\sqrt{\frac{2(n-1)}{\varepsilon}}\nonumber\\
k<0,&& \text{for }|z|<\sqrt{\frac{2(n-1)}{\varepsilon}}\,.\nonumber
\end{eqnarray}
\end{remark}

\subsection{Starshapedness}~

A curve $z:S^1\to\complex{*}$
will be called {\sl starshaped}, if $\langle\nu,z\rangle>0$ everywhere. In this
section we want to show that starshapedness is preserved as long as $r=|z|>0$.
To this end we need some evolution equations and the following lemma which
will be used several times in our computations

\begin{lemma}\label{lemma  10}
On any curve $\gamma\in\complex{*}$ we have
\begin{eqnarray}
1&=&|\nabla r|^2+\langle\nu,e_r\rangle^2,\label{eq 10a}\\
\Delta r&=&\langle\nu,e_r\rangle
\left(\frac{\langle\nu,e_r\rangle}{r}-k\right),\label{eq 10b}\\
\nabla\langle\nu,e_r\rangle&=&\left(k-\frac{\langle\nu,e_r\rangle}{r}\right)
\nabla r.\label{eq 10c}
\end{eqnarray}
\end{lemma}
\begin{proof}
We have 
$$z=\frac{\langle z,z'\rangle}{|z'|^2} z'+\langle z,\nu\rangle\nu$$
and therefore
$$r^2=|z|^2=\frac{|\langle z,z'\rangle|^2}{|z'|^2}+\langle z,\nu\rangle^2=
\left|\nabla\frac{r^2}{2}\right|^2+\langle z,\nu\rangle^2$$
which proves the first equation. The second equation then follows
from the first and Gau\ss' equation $\Delta z=-k\nu$. The last equation 
follows from $\nu'=kz'$.
\end{proof}

Next we want to compute the evolution equations for $r$ and 
$\langle\nu, e_r\rangle$.
\begin{lemma}\label{lemma 5}
The following equations are valid
\begin{equation}
\dt r=\Delta r+\frac{n-1}{r}\langle\nabla r,e_r\rangle
-\frac{1}{r}\langle\nu,e_r\rangle^2
-(n-1)\frac{1}{r}\label{evol 5b}
\end{equation}
\begin{eqnarray}
\dt\nur
&=&\Delta\nur+\frac{n-1}{r}\langle\nabla\nur, e_r\rangle\label{evol 7b}\\
&&\hspace{-20pt}-2n\frac{\nur}{r^2}(1-\nur^2)+\left(k-\frac{\nur}{r}\right)^2
\nur\nonumber
\end{eqnarray}
\end{lemma}
\begin{proof}
We first observe that
$$\dt r=\frac{1}{r}\langle z,\dt z\rangle=-f\langle\nu,
e_r\rangle,$$
so that in view of (\ref{eq 10a}) and 
$\langle\nabla r,e_r\rangle=|\nabla r|^2$
\begin{equation}\nonumber
\dt r=-k\langle\nu,e_r\rangle+\frac{n-1}{r}\langle\nabla r,e_r\rangle
-\frac{n-1}{r}.
\end{equation}
On the other hand, by (\ref{eq 10b}) the Laplacian of $r$ is given by
$$\Delta r=\langle\nu,e_r\rangle
\left(\frac{\langle\nu,e_r\rangle}{r}-k\right).$$
Hence
\begin{equation}
\dt r=\Delta r+\frac{n-1}{r}\langle\nabla r,e_r\rangle
-\frac{1}{r}\langle\nu,e_r\rangle^2
-(n-1)\frac{1}{r}\,.\nonumber
\end{equation}
{F}or the evolution equation of $\nur$ we compute
\begin{eqnarray}
\dt\nur&=&\langle\nabla f,e_r\rangle-\frac{\nur}{r^2}\langle z,\dt z\rangle
-\frac{f}{r}\nonumber\\
&=&\langle\nabla f,e_r\rangle-\frac{f}{r}(1-\nur^2)\nonumber\\
&=&\langle\nabla k,e_r\rangle+\frac{n-1}{r}\langle\nabla\nur, e_r\rangle
\nonumber\\
&~&-\left(\frac{f}{r}+(n-1)\frac{\nur}{r^2}\right)(1-\nur^2).\nonumber
\end{eqnarray}
To compute the Laplacian we use (\ref{eq 10b}) and (\ref{eq 10c}) to
derive
\begin{eqnarray}
\Delta\nur&=&\langle\nabla k,e_r\rangle-\left(k-\frac{\nur}{r}\right)^2\nur
\nonumber\\
&~&-\frac{1-\nur^2}{r}\left(k-\frac{\nur}{r}\right)+\frac{\nur}{r^2}(1-\nur^2).
\nonumber
\end{eqnarray}
Combining these two equations we get
\begin{eqnarray}
\dt\nur
&=&\Delta\nur+\frac{n-1}{r}\langle\nabla\nur, e_r\rangle\nonumber\\
&~&-\left(\frac{f}{r}+n\frac{\nur}{r^2}\right)(1-\nur^2)\nonumber\\
&~&+\left(k-\frac{\nur}{r}\right)^2\nur+
\frac{1-\nur^2}{r}\left(k-\frac{\nur}{r}\right)\nonumber\\
&=&\Delta\nur+\frac{n-1}{r}\langle\nabla\nur, e_r\rangle\nonumber\\
&~&-2n\frac{\nur}{r^2}(1-\nur^2)+\left(k-\frac{\nur}{r}\right)^2\nur.\nonumber
\end{eqnarray}
\end{proof}

\begin{corollary}\label{coral 3}
If the initial curve $\gamma_0\subset\complex{*}$
is starshaped, i.e. if 
$$\nur>0$$
for all $p\in\gamma_0$, then this remains true as long as a smooth solution
of (\ref{ecsf}) exists in $\complex{*}$.
\end{corollary}
\begin{proof}
Suppose a smooth solution in $\complex{*}$ exists on a time interval
$[0,t_1)$ and that $0<t_0<t_1$ is some fixed time.
There exists a constant $\lambda>0$ depending on $t_0$ such that
$$\frac{2n}{r^2}(1-\nur^2)-\left(k-\frac{\nur}{r}\right)^2\le\lambda$$
holds on the compact time interval $[0,t_0]$. Then by (\ref{evol 7b})
we can estimate for $t\in [0,t_0]$
$$\dt\nur\ge\Delta\nur+\frac{n-1}{r}\langle\nabla\nur,e_r\rangle-\lambda\nur$$
and the maximum principle implies that
$$\nur\ge\left(\inf_{\gamma_0}\nur\right)e^{-\lambda t}\,,\quad
\forall\,t\in[0,t_0].$$
\end{proof}
\begin{remark}\label{rem b}~

\begin{itemize}
\item[(a)]
Corollary \ref{coral 3} is sharp in the following sense: 
If $n=1$, i.e. if we consider the usual curve shortening flow, and if the
initial curve is a circle of some radius including the origin but not 
centered at the origin, then this remains starshaped as long as the origin 
is contained in the shrinking circles. As soon as one circle touches the 
origin, starshapedness will fail.
\item[(b)]
If $n>1$, then by Lemma \ref{lemma long} a smooth solution of (\ref{ecsf}) 
with initial curve $\gamma_0\in\complex{*}$ is contained in $\complex{*}$ 
for any $t\in[0,T_{\text{sing}})$.
\end{itemize}
\end{remark}

We can prove even stronger estimates:
\begin{lemma}\label{lemma 3b}
If the initial curve $\gamma_0\subset\complex{*}$
is starshaped, then for any $p\ge 2$ there exists a constant $\epsilon_p>0$ 
such that 
\begin{equation}\nonumber
r^{-p}\nur\ge\epsilon_p
\end{equation}
holds true as long as a smooth solution of (\ref{ecsf}) exists in 
$\complex{*}$.
\end{lemma}

\begin{proof}
{F}irst we compute the evolution equation of $m:=r^{-p}\nur$.
Therefore, in general the following computations
are valid only as long as $r>0$. According to Corollary \ref{coral 3} and 
Remark \ref{rem b}, this quantity is well defined for $t\in[0,T_{\text{sing}})$, if 
$n>1$ and in case $n=1$ as long as $r>0$. Lemma \ref{lemma 5} implies
\begin{eqnarray}
\dt m
&=&\nur
\left\{\lap r^{-p}+\frac{n-1}{r}\langle\nabla r^{-p},e_r\rangle
+p(p+2)r^{-p-2}\nur^2\right.\nonumber\\
&&\left.-p(p-n+2)r^{-p-2}\right\}
+r^{-p}
\left\{\Delta\nur\phantom{\frac{1}{r}}
\right.\nonumber\\
&&\left.+\frac{n-1}{r}\langle\nabla\nur, e_r\rangle
-2n\frac{\nur}{r^2}(1-\nur^2)\right.\nonumber\\
&&\left.+\left(k-\frac{\nur}{r}\right)^2\nur\right\}.\nonumber
\end{eqnarray}
The Laplacian of $m$ is given by
\begin{eqnarray}
\lap m
&=&r^{-p}\Delta\nur+\nur\Delta r^{-p}\nonumber
+2\langle\nabla\nur,\nabla r^{-p}\rangle\nonumber\\
&=&r^{-p}\Delta\nur+\nur\Delta r^{-p}\nonumber\\
&&-2pr^{-p-1}\left(k-\frac{\nur}{r}\right)(1-\nur^2),\nonumber
\end{eqnarray}
where we have used (\ref{eq 10a}) and (\ref{eq 10c}).
Therefore it follows
\begin{eqnarray}
\dt m
&=&\lap m+\frac{n-1}{r}\langle\nabla m,e_r\rangle
+2pr^{-p-1}\left(k-\frac{\nur}{r}\right)(1-\nur^2)\nonumber\\
&&+m\left\{-\,\frac{p^2+2p+2n}{r^2}\,
(1-\nur^2)+\frac{np}{r^2}+\left(k-\frac{\nur}{r}\right)^2\right\}.\nonumber
\end{eqnarray}
We will show that $m$ cannot admit a decreasing positive minimum.
At a minimum of $m$ we have $\nabla m=0,\Delta m\ge 0$. From
$$\nabla m=r^{-p}\left(k-(p+1)\frac{\nur}{r}\right)\nabla r=0$$
we conclude, that at such a point either $\nabla r=0$, or 
$k-\frac{\nur}{r}=p\frac{\nur}{r}$.
In the first case we also have $1-\nur^2=0$ and then
$$\dt m\ge m\left\{\frac{np}{r^2}+\left(k-\frac{\nur}{r}\right)^2\right\}>0.$$
In the second case, we substitute $k-\frac{\nur}{r}$ and obtain
\begin{eqnarray}
\dt m&\ge& m\left\{\frac{p^2-2p-2n}{r^2}(1-\nur^2)
+\frac{np}{r^2}+p^2\frac{\nur^2}{r^2}\right\}\nonumber\\
&=&m\left\{2(p+n)\frac{\nur^2}{r^2}
+\frac{(p-2)(p+n)}{r^2}\right\}\ge 0\nonumber
\end{eqnarray}
since $p\ge 2$. This proves the lemma.
\end{proof}

We can also prove, that a weaker form of starshapedness is preserved.

\begin{definition}
A curve $\gamma\in\complex{*}$ will be called {\sl austere}, if
$$\nur>-1$$
for all points on $\gamma$.
\end{definition}

This means that the angle $\beta$ defined by $\cos\beta=\nur$
is contained in the open interval $(-\pi,\pi)$.

\begin{lemma}\label{lemma austere}
If the initial curve $\gamma_0\subset\complex{*}$
is austere, then this holds true as long as a smooth solution 
of (\ref{ecsf}) exists in $\complex{*}$.
\end{lemma}

\begin{proof}
{F}rom equations (\ref{eq 10a}), (\ref{eq 10c}) and the evolution equation of
$\nur$ one easily deduces the evolution equation of $\beta$, i.e.
\begin{equation}\label{eq beta}
\dt\beta=\Delta\beta+\frac{n-1}{r}\langle\nabla\beta,e_r\rangle
+\frac{n}{r^2}\sin(2\beta).
\end{equation}
Then, using the periodicity of $\sin(2\beta)$ and the estimate $\sin x\le \frac{2x}{\pi}$, for $-\frac{\pi}{2}\le x\le 0$
we conclude 
$$\dt(\beta-\pi)\le\Delta(\beta-\pi)
+\frac{n-1}{r}\langle\nabla(\beta-\pi),e_r\rangle
+\frac{4n}{\pi r^2}(\beta-\pi)$$
for all $-\frac{\pi}{4}\le\beta-\pi\le0$.
The parabolic maximum principle implies that $\beta-\pi$ stays negative as
long as $\frac{n}{r^2}$ stays bounded. The estimate $\beta>-\pi$ can be
derived similarly.
\end{proof}

In view of Lemma \ref{lemma austere} we also add the following lemma.

\begin{lemma}
If $z:S^1\to\complex{*}$ is a closed, austere and tamed curve, then $z$ is
starshaped and the function
$$\delta:=r^n\nur$$
satisfies
\begin{equation}\label{comparison}
\delta\ge\rho^n,
\end{equation}
where
$$\rho:=\min_{\phi\in S^1}r(\phi).$$
\end{lemma}

\begin{proof}
By (\ref{eq 10c}) we have
\begin{equation}\nonumber
d\delta=r^nfdr
\end{equation}
and since $f>0$ we must have $|\nabla r|^2=1-\nur^2=0$ at a minimum of 
$\delta$.
Since $z$ is austere, we conclude that at such a point $\nur=1$ and that the
minimum of $\delta$ is bounded from below by $\rho^n>0$. Hence $z$ must
also be starshaped.
\end{proof}

\section{Singularities for starshaped profile curves}\label{sing}
In this section we will prove Theorems \ref{main 1} and \ref{main 3}.
\pproof{\ref{main 1}}
{F}irst we recall that $r>0$ on $[0,T_{\text{sing}})$ since $n\ge 2$ (see Lemma \ref{lemma long}).
We will prove that $\sup|f|<\infty$ as long as $\inf r>0$. In view of 
Lemma \ref{lemma long} this will prove the theorem.
Let $\beta$ be the angle determined by $\cos\beta=\nur$
and suppose $p=p(\beta)$ is an arbitrary smooth function
of $\beta$ which will be determined later. We will denote a partial derivative
w.r.t. $\beta$ by a prime.
The evolution equation for $s:=fp$ is given by
\begin{eqnarray}
\dt s
&=&\Delta s+\frac{n-1}{r}\langle\nabla s,e_r\rangle
-2\langle\nabla p,\nabla f\rangle\nonumber\\
&&+s\left(k^2-\frac{p''}{p}\left(k-\frac{\nur}{r}\right)^2
\right.\nonumber\\
&&\left.+\frac{n-1}{r^2}(2\nur^2-1)+\frac{np'}{r^2p}\sin(2\beta)\right)
\nonumber\\
&=&\Delta s+\frac{n-1}{r}\langle\nabla s,e_r\rangle
-\frac{2}{p}\langle\nabla p,\nabla s\rangle\nonumber\\
&&+s\left(k^2+\left(2\left(\frac{p'}{p}\right)^2-\frac{p''}{p}\right)
\left(k-\frac{\nur}{r}\right)^2
\right.\nonumber\\
&&\left.+\frac{n-1}{r^2}(2\nur^2-1)+\frac{np'}{r^2p}\sin(2\beta)\right)
\nonumber.
\end{eqnarray}
Suppose $\liminf_{t\to T_{\text{sing}}}\left(\min_{\gamma_t}r\right)>0$. Then by
Lemma \ref{lemma 3b} we can find $\sigma >1, \epsilon>0$ such that 
$\cos(\sigma\beta)\ge\epsilon>0$ holds on $[0,T_{\text{sing}})$. 
Now we choose $p=\frac{1}{\cos(\sigma\beta)}$. 
Then $2\left(\frac{p'}{p}\right)^2-\frac{p''}{p}=-\sigma^2$
gives
\begin{eqnarray}
\dt s
&=&\Delta s+\frac{n-1}{r}\langle\nabla s,e_r\rangle
-\frac{2}{p}\langle\nabla p,\nabla s\rangle\nonumber\\
&&+s\left(k^2-\sigma^2
\left(k-\frac{\nur}{r}\right)^2
\right.\nonumber\\
&&\left.+\frac{n-1}{r^2}(2\nur^2-1)+\frac{n\sigma}{r^2}\tan(\sigma\beta)
\sin(2\beta)\right)
\nonumber.
\end{eqnarray}
Again, since $r$ is bounded from below and since $\sigma>1$, 
there exists a small positive constant $\rho$ and a positive constant 
$c$ depending on $n, \sigma, \epsilon, \inf r$ such that 
\begin{eqnarray}
&&k^2-\sigma^2
\left(k-\frac{\nur}{r}\right)^2\nonumber\\
&&+\frac{n-1}{r^2}(2\nur^2-1)+\frac{n\sigma}{r^2}\tan(\sigma\beta)
\sin(2\beta)\nonumber\\
&&\le-\rho s^2+c.\nonumber
\end{eqnarray}
Substituting this in the evolution equation for $s$, 
we may apply the maximum principle on $[0,T_{\text{sing}})$ 
to see that $s$ must be bounded from above and from below on the 
time interval $[0,T_{\text{sing}})$. 
This implies that $|f|$ and $|k|$ are 
uniformly bounded on $[0,T_{\text{sing}})$ as well
which contradicts $T_{\text{sing}}<\infty$.
\endpproof

Now let us assume that $\gamma\in\complex{}^*$ is a smooth, starshapedly 
immersed curve with winding number $\omega_0:=\wind_0(z)$. 
$\gamma$ can be parameterized by a smooth map
\begin{eqnarray}
z&:&S^1\to\complex{}\nonumber\\
z(\phi)&:=&r(\phi)\binom{\cos(\omega_0\phi)}{\sin(\omega_0\phi)},\nonumber
\end{eqnarray}
where $r$ is a positive
smooth function on $S^1\cong[0,2\pi)$.
We abbreviate a derivative w.r.t. $\phi$ by a prime. Then 
\begin{equation}\label{curveeq 1}
z'(\phi)=r'(\phi)\binom{\cos(\omega_0\phi)}{\sin{(\omega_0\phi)}}
+\omega_0r(\phi)\binom{-\sin{(\omega_0\phi)}}{\cos{(\omega_0\phi)}},
\end{equation}
\begin{equation}\label{curveeq 2}
\nu(\phi)=\frac{1}{\sqrt{g(\phi)}}
\left(\omega_0r(\phi)\binom{\cos(\omega_0\phi)}{\sin{(\omega_0\phi)}}
-r'(\phi)\binom{-\sin{(\omega_0\phi)}}{\cos{(\omega_0\phi)}}
\right),
\end{equation}
where $\nu$ is the outward pointing unit normal and $g$ denotes the induced
metric on the curve given by
\begin{equation}\label{curveeq 3}
g(\phi)=\langle z'(\phi),z'(\phi)\rangle=(\omega_0r(\phi))^2+(r'(\phi))^2.
\end{equation}
The formula for the curvature $k$ is
\begin{equation}\label{curveeq 4}
k=\frac{1}{g}\langle z'',\nu\rangle=\frac{\omega_0}{\sqrt{g}}\left(1+
\frac{(r')^2-rr''}{g}\right).
\end{equation}
Moreover,
\begin{equation}\label{curveeq 5}
\frac{\langle z,\nu\rangle}{|z|^2}=\frac{\omega_0}{\sqrt{g}}
\end{equation}
and 
\begin{equation}\label{curveeq 6}
f=k+(n-1)\frac{\langle z,\nu\rangle}{|z|^2}=\frac{\omega_0}{\sqrt{g}}
\left(n+\frac{(r')^2-rr''}{g}\right).
\end{equation}
Let us define the angle $\beta$ by
\begin{equation}\label{beta}
\beta:=\frac{1}{\omega_0}\arctan{\frac{r'}{\omega_0r}}.
\end{equation}
We obtain
\begin{equation}\label{beta2}
\beta'=\frac{rr''-(r')^2}{(\omega_0r)^2+(r')^2}
\end{equation}
so that with (\ref{curveeq 6}) we get
\begin{equation}\label{curveeq 11}
f\sqrt{g}=\omega_0(n-\beta').
\end{equation}
$\beta$ measures how much the curves differ from a circle centered at the
origin. If $\gamma$ is tamed, then $f>0$. This is equivalent to  
$\left(n\phi-\beta\right)'>0$.

\pproof{\ref{main 3}}
\begin{itemize}
\item[(i)]
The stability of $\sE$ follows directly from Proposition \ref{lemma mon},
Lemma \ref{lemma emb}, Corollaries \ref{coral 1}, \ref{coral 2}
and \ref{coral 3}. The existence of type-1 singularities within $\sE$
will be shown below and the existence of type-2 singularities follows from
Neves' work \cite{neves} in the following way:
As pointed out in \cite{neves}, a straight line in $\complex{*}$ 
will develop a type-2 singularity in case $n\ge 2$ (cf. Figure \ref{fig 0}). 
\begin{figure}[h]
\begin{center}
\picture{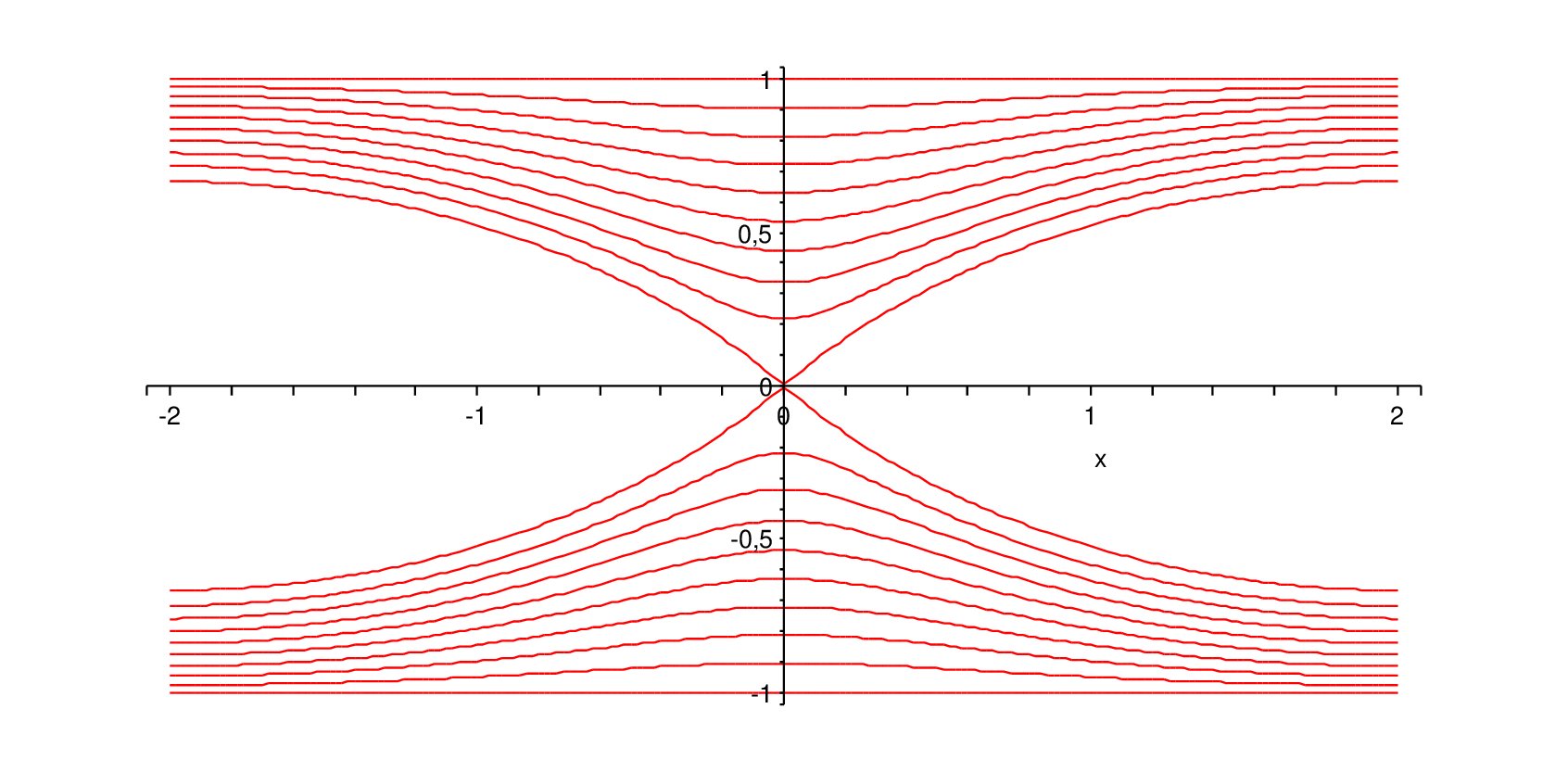}{1}
\caption{\sl Straight lines in $\complex{*}$ develop a type-2 singularityin case $n\ge 2$.}
\label{fig 0}
\end{center}
\end{figure}
Thus, by some glueing technique,
we obtain also examples of closed profile curves, that generate Lagrangian 
submanifolds in $\sE$ and which will develop a type-2
singularity at the origin as $t\to T_{\text{sing}}$ and for which 
$T_{\text{sing}}<\frac{1}{\varepsilon}$ so that the enclosed area 
does not vanish as $t\to T_{\text{sing}}$.
{F}igure \ref{fig 1} depicts such a singularity formation. 
\begin{figure}[h]
\begin{center}
\picture{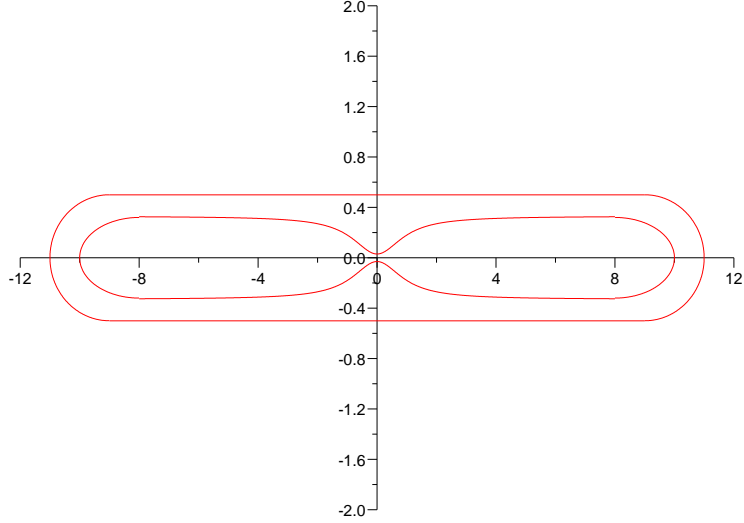}{1}
\caption{\sl The outer curve will evolve into a curve that roughly
looks like the inner one. A type-2 singularity forms.}
\label{fig 1}
\end{center}
\end{figure}

\item[(ii)]
Let $\gamma$ be an arbitrary profile curve such that the corresponding
Lagrangian submanifold $L$ belongs to $\sE$ and such that $\gamma$
is symmetric by rotations of ${2\pi}/{l}$ with $l\ge 1+4n\omega_0$. 
Let $\phi_0\in[0,2\pi)$ be a point with
$r'(\phi_0)=0$. W.l.o.g. we may assume $\phi_0=0$. Then by symmetry we must
have $r'({2\pi}/{l})=0$ as well. Since $f>0$ the function $n\phi-\beta$
is strictly increasing and from $\beta(0)=\beta({2\pi}/{l})=0$ we conclude
$$0\le n\phi-\beta(\phi)\le \frac{2n\pi}{l}$$
for all $\phi\in[0,2\pi/l)$. Hence
$\beta(\phi)\le\frac{2n\pi}{l}<\frac{\pi}{2\omega_0}-\frac{\pi}{2\omega_0(1+4n\omega_0)}$ and
$\beta(\phi)\ge n\phi-\frac{2n\pi}{l}>\frac{\pi}{2\omega_0(1+4n\omega_0)}-\frac{\pi}{2\omega_0}$.
This implies that there exists a constant $c>0$ independent of $\gamma$
such that
$$c-\frac{\pi}{2}\le\arctan\frac{r'}{\omega_0 r}\le \frac{\pi}{2}-c$$
and therefore 
\begin{equation}\label{gradient}
\sup_{\gamma}\left|\frac{r'}{r}\right|\le C
\end{equation}
for some uniform constant $C$ independent of $\gamma$.  
Let $r_+(\gamma):=\max_{\gamma}r$, $r_-(\gamma):=\min_{\gamma}r$. Integrating
(\ref{gradient}) on $\gamma$ we obtain the Harnack estimate
\begin{equation}\label{harnack}
\frac{r_+(\gamma)}{r_-(\gamma)}\le e^{2\pi C},
\end{equation}
which is independent of $\gamma$.\\

\noindent
Now for some $\gamma_0$ as above
let $\gamma_t$ be the profile curves given by the equivariant mean 
curvature flow (\ref{ecsf}). Since the mean curvature flow preserves 
the invariance under isometries we conclude from the first part of the 
proof that $\gamma_t$ stays in our class and consequently
(\ref{harnack}) holds for all $t\in[0,T_{\text{sing}})$.\\

\noindent
{F}rom Theorem \ref{main 1} we get $\liminf_{t\to T_{\text{sing}}}r_-(t)=0$
and (\ref{harnack}) then implies
$$\limsup_{t\to T_{\text{sing}}}r_+(t)=0.$$
In particular, $\lim_{t\to T_{\text{sing}}} A(\gamma_t)=0$ which in view
of the symplectic area formula and the monotonicity of $\gamma_0$
is only possible, if $T_{\text{sing}}=\frac{1}{\varepsilon}$.\\

\noindent
By Proposition \ref{prop main 1} the rescaled flow 
$$\tilde z=\frac{z}{\sqrt{2(T_{\text{sing}}-t)}}$$
becomes Hamiltonian and the symplectic area 
$$\tilde A(\tilde z)=\frac{1}{2}\int\limits_{S^1}\langle\tilde z,\tilde\nu
\rangle d\tilde\mu$$ is fixed which is equivalent to
$$\frac{\omega_0}{2}\int\limits_{0}^{2\pi}\tilde r^2$$ being fixed. 
By the mean value theorem and the logarithmic Harnack inequality
(\ref{harnack}) for $r$ and $\tilde r$
we then conclude that there exist uniform positive upper and lower bounds for 
$\tilde r$. Similar as in the proof of Theorem \ref{main 1} we can then
derive a uniform upper bound for $|\tilde f|$. From the definition of
the rescaling it follows that $|\tilde f|$ is uniformly bounded, if and 
only if $|A|^2\le\frac{c}{T_{\text{sing}}-t}$ for some constant $c$,
so that the singularity must be of type-1. From Huisken's result
for the blow-up limits of type-1 singularities we then obtain a smooth 
closed limiting curve $\tilde\gamma_\infty\subset\complex{*}$
that is one of the self-similarly contracting solutions classified
by Anciaux.
\end{itemize}
\endpproof

\bibliographystyle{amsplain}

\end{document}